\documentclass[12pt]{amsart}       
\usepackage{amssymb, amscd}

\newtheorem{theorem}{Theorem}[section]
\newtheorem{lemma}[theorem]{Lemma}
\newtheorem{proposition}{Proposition}[section] 
\newtheorem{corollary}[theorem]{Corollary}  

\theoremstyle{definition}

\theoremstyle{remark}
\newtheorem{remark}[theorem]{Remark}

\numberwithin{equation}{section}

\renewcommand{\a}{\alpha}
\newcommand{\ch}{\mbox{ch}}
\newcommand{\Cs}{\mathbb C^{\times}}

\newcommand{\FGC}{{\mathcal F}_{\G\times \mathbb C^{\times}} }

\newcommand{\FGGC}{\overline{\mathcal F}_{\G\times \mathbb C^{\times}}}
\newcommand{\tFGGC}{\widetilde{\mathcal F}_{\G\times \mathbb C^{\times}}} 

\newcommand{\G}{\Gamma}
\newcommand{\GC}{\Gamma\times \mathbb C^{\times}}

\newcommand{\Gm}{{\Gamma}_m}
\newcommand{\GCm}{{\Gamma}_m\times\mathbb C^{\times}}
\newcommand{\Gn}{{\Gamma}_n}
\newcommand{\GCn}{{\Gamma}_n\times\mathbb C^{\times}}
\newcommand{\g}{\gamma}
\newcommand{\hg}{\widehat{\mathfrak h}_{\G, \wt}}
\newcommand{\la}{\lambda}
\newcommand{\loopg}{\widehat{\mathfrak g}}
\newcommand{\qaa}{U_q(\widehat{\mathfrak g})}
\newcommand{\qta}{U_q(\widehat{\widehat{\mathfrak g}})}
\newcommand{\ta}{\widehat{\widehat{\mathfrak g}}}
\newcommand{\qtap}{U_{q, p}(\widehat{\widehat{\mathfrak g}})}
\newcommand{\LG}{ {\Lambda}_{\G}}

\newcommand{\ep}{\epsilon}
\newcommand{\RC}{R(\mathbb C^{\times})}
\newcommand{\RG}{R_{\G}}

\newcommand{\RGGC}{\overline{R}_{ \G\times \mathbb C^{\times}}}

\newcommand{\RGC}{R_{\G\times \mathbb C^{\times}}}

\newcommand{\Rz}{R_{\mathbb Z}(\Gamma)}
\newcommand{\Rzz}{\overline{R}_{\mathbb Z}({\G})}
\newcommand{\s}{\sigma}
\newcommand{\SG}{ S_\G }
\newcommand{\SGC}{S_{\GC}}

\newcommand{\SGGC}{\overline{S}_{\G\times\mathbb C^{\times} }}

\newcommand{\VG}{V_{ \G}}
\newcommand{\VGC}{V_{\G\times \mathbb C^{\times}}}

\newcommand{\tVGGC}{\widetilde{V}_{ \G\times\mathbb C^{\times}}}

\newcommand{\wt}{\xi}

\begin{document}

\title[Quantum vertex operators and McKay correspondence]
    {Quantum vertex representations via finite groups and 
the McKay correspondence}
\author{Igor B. Frenkel}
\address{Frenkel: Department of Mathematics,
   Yale University, 
   New Haven, CT 06520}
\author{Naihuan Jing}
\address{Jing: Department of Mathematics,
   North Carolina State Univer\-sity,
   Ra\-leigh, NC 27695-8205 
\qquad
Mathematical Sciences Research Institute, 1000 Centennial Drive,
Berkeley, CA 94720
}
\email{jing@math.ncsu.edu}
\thanks{I.F. is supported in part by NSF grant 
DMS-9700765. N.J. is supported in part by NSA grant 
MDA904-97-1-0062 and NSF grant DMS-9970493.}
\author{Weiqiang Wang}
\address{Wang: Department of Mathematics,
   North Carolina State Univer\-sity,
   Ra\-leigh, NC 27695-8205 
\qquad
Department of Mathematics,
   Yale University,
   New Haven, CT 06520}
\email{wqwang@math.ncsu.edu}


\begin{abstract}
We establish a $q$-analog of our recent work on
vertex representations and the McKay correspondence.
For each finite group $\G$ we construct a Fock space
and associated vertex operators  
in terms of wreath products of $\Gamma\times \mathbb C^{\times}$
and the
symmetric groups.
 An important special
case is obtained when $\Gamma$ is a finite subgroup
of $SU_2$, where our construction yields a group theoretic
realization of the representations of
the quantum affine and quantum toroidal algebras
of $ADE$ type.

\end{abstract}

\maketitle

\section{Introduction} \label{S:intro}

In our previous paper \cite{FJW} (see \cite{W, FJW} for historical
remarks and motivations)  
we have shown that the
 basic representation
of an affine Lie algebra $\loopg$ of ADE type
can be constructed from a finite subgroup $\G$ of $SU_2$ related to the 
Dynkin diagram of $\loopg$ via the McKay correspondence.
In particular, we have recovered a well-known
construction \cite{FK, Se} of the basic representation of $\loopg$
from the root lattice $Q$ of the corresponding finite dimensional
Lie algebra $\mathfrak g$. In fact our construction yields
naturally the vertex representation of the toroidal
Lie algebra $\ta$ which contains the affine Lie algebra as
a distinguished subalgebra.

The main goal of the present paper is to $q$-deform our construction 
in \cite{FJW}.
Again as in the undeformed case we will naturally obtain
the earlier construction \cite{FJ} of the basic representation
of the quantum affine algebra 
$\qaa$ from the root lattice $Q$ and its generalization
to the quantum toroidal algebra $\qta$ \cite{GKV} (also
cf. \cite{Sa, J3}). The $q$-deformation is achieved
by replacing consistently the  representation theory of $\G$ by that
of $\GC$. The representation ring for $\Cs$ is identified with the ring of
Laurent polynomials $\mathbb C[q, q^{-1}]$ so that the formal variable
$q$ corresponds to the natural one-dimensional representation
of $\Cs$. It turns out that rather complicated expressions for operators
in Drinfeld realization of the quantum affine algebra $\qaa$ and the 
quantum toroidal algebra $\qta$ follow 
instantly from the simple extra factor $\Cs$. The idea to use representations of $\Cs$
to obtain a $q$-deformation of the basic representation was mentioned
in \cite{Gr} and is widely used in geometric constructions of 
representations (see e.g. \cite{CG}).

As in the previous paper \cite{FJW} we give the construction of 
quantum vertex operators starting from an arbitrary finite group
$\G$ and a self-dual virtual character $\wt$ of $\GC$. Using
the restriction and induction functors in representation theory
of wreath products of $\GC$ with the symmetric group
$S_n$ for all $n$ we construct two ``halves'' of quantum vertex operators
corresponding to any irreducible character $\g$ of $\GC$. Then choosing an 
irreducible character of $\Cs$, i.e. an integer power of $q$ we
assemble both halves into one quantum vertex operator. 

The special case when $\G$ is a subgroup of $SU_2$ is important for the
application to representation theory of $\qaa$ and for
relations \cite{W} to the theory of Hilbert schemes of points
on surfaces.
 To recover the basic representation of 
$\qaa$ we choose 
$$\wt=\g_0\otimes (q+q^{-1})-\pi\otimes 1_{\Cs},
$$
where $\g_0$ and $1_{\Cs}$ are the trivial characters
of $\G$ and $\Cs$ respectively, $q$ and $q^{-1}$ are
 the natural and its dual 
characters of $\Cs$, and $\pi$ is the natural character
of $\G$ in $SU_2$. 
The fact that the quantum toroidal algebra intrinsically presents
in our construction is an additional indication of its importance
in representation theory of quantum affine algebras. 
Moreover when $\G$ is cyclic of order $r+1$, 
$\pi\simeq \g\oplus \g^{-1}$, where $\g$ is
the natural character of $\G$, 
one can modify our virtual character $\wt$
with an extra parameter $p=q^k, k\in\mathbb Z$ by letting  
$$
\wt=\g_0\otimes (q +q^{-1})-
(\g\otimes p+\g^{-1}\otimes p^{-1}). 
$$
In the special case when $p=q^{\pm 1}$ the quantum vertex representation
of the quantum toroidal algebra $\qta$ can be factored
to the basic representation of the quantum affine algebra $\qaa$.
This is a $q$-analog of the factorization in the undeformed case, which exists
for an arbitrary simply-laced affine Lie algebra.

To obtain the basic representations of quantum toroidal and affine
algebras we only need the quantum vertex operators corresponding to
irreducible representations of $\G$ and their negatives in the
Grothendieck ring of this group. We attach
two halves  of quantum vertex operators using the simplest nontrivial
representations of $\Cs$ namely $q$ and $q^{-1}$. Each
of the two choices and only these two yield the basic representations
of $\qta$ and $\qaa$, in a perfect correspondence with the
construction in \cite{FJ}. This choice of an irreducible character
of $\Cs$ is essentially the only freedom that exists in our construction
of quantum vertex operators for the quantum affine and toroidal algebras 
 and is fixed by
comparison with the algebra relations. However it raises the question
of constructing a ``natural'' quantum vertex operator corresponding to any 
virtual character  $\g$ of $\G$. This question is closely related to
the well-known problem of finding a $q$-deformation of 
vertex operator algebras associated to the basic representation of an affine
Lie algebra. 

This paper is organized in a way similar to \cite{FJW}.
In Sect. \ref{sect_wreath} we review the theory of wreath products of
$\G$ and extend it to $\GC$. In Sect. \ref{sect_weight} we define
the weighted bilinear form on $\GC$ and its wreath products. In Sect. \ref{sect_mckay} we introduce two 
distinguished $q$-deformed weight functions associated to subgroups of $SU_2$.
In Sect. \ref{sect_heis} we define the Heisenberg
algebra associated to $\G$ and the weighted bilinear form, and we
construct its representation in a Fock space. In Sect. \ref{sect_isom}
we establish the isometry between the representation ring
of wreath products of $\GC$ and the Fock space
representation of the Heisenberg algebra. In Sect. \ref{sect_vertex}
we construct quantum vertex operators acting on the representation
ring of the wreath products.
In Sect. \ref{sect_ade} we 
obtain the basic representations of quantum toroidal algebras and
quantum affine algebras
from representation theory of wreath products for $\GC$.

\section{Wreath products and vertex representations} \label{sect_wreath}

\subsection{The wreath product $\Gn$} 
  Let $\Gamma$ be a finite group and $n$ a non-negative integer.
The wreath product $\Gn$ is the semidirect product of the $n$-th direct product $\G^n=\G\times\cdots
\times\G$ and the symmetric group $S_n$:
$$
 \Gamma_n = \{(g, \sigma) | g=(g_1, \ldots, g_n)\in {\Gamma}^n,
\sigma\in S_n \}
$$
 with the group multiplication
$$
(g, \sigma)\cdot (h, \tau)=(g \, {\sigma} (h), \sigma \tau ) ,
$$
where $S_n$ acts on $\G^n$ by permuting the factors.

Let $\G_*$ be the set of conjugacy classes of $\G$ consisting of 
$c^0=\{1\}$, $c^1$, $\dots$,
$c^r$ and $\G^*$ be the set of 
$r+1$ irreducible characters: $\g_0, \g_1, \dots, \g_r$. Here we denote
the trivial character of $\G$ by $\g_0$. The order of the centralizer
of an element in the conjugacy class
$c$ is denoted by $\zeta_c$, so the order of the conjugacy class $c$
is $|c|=|\G|/\zeta_c$, where $|\G|$ is the order of $\G$.

 A partition $\la=(\la_1, \la_2, \ldots, \la_l)$ 
is a decomposition of $n=|\la|=\la_1+\cdots+\la_l$ with nonnegative integers:
$\la_1\geq \dots \geq \la_l \geq 1$, where
$l=l (\la )$ is called the {\em length} of the partition
$\la $ and $\la_i$ are called the {\em parts} of $\la$. 
Another notation for $\la$ is 
$$
\la=(1^{m_1}2^{m_2}\cdots)
$$
with $m_i$ being the multiplicity of parts equal to $i$ in 
$\la$. 
Denote by $\mathcal P$ the set of all partitions of integers and by $\mathcal P(S)$ the set of
all partition-valued functions on a set $S$.
The weight of a partition-valued function $\rho=(\rho(s))_{s\in S}$ is defined to be $\|\rho\|=\sum_{s\in S}|\rho(s)|$. We also denote
by $\mathcal P_n$ (resp. $\mathcal P_n(S)$) the subset
of $\mathcal P$ (resp. $\mathcal P(S)$) of partitions
with weight $n$.

Just as the conjugacy classes of $S_n$ are parameterized by
partitions,
the conjugacy classes of $\Gn$ are parameterized by partition-valued
functions on $\G_*$. Let $x=(g, \s)\in {\Gamma}_n$, where
$g=(g_1, \ldots, g_n) \in {\Gamma}^n$ and $\s\in S_n$ is
presented as a product of disjoint cycles. For each 
cycle $(i_1 i_2 \cdots i_k)$ of $\sigma$, we define the {\em cycle-product}
element $g_{i_k} g_{i_{k -1}}
\cdots g_{i_1} \in \Gamma$, which is determined up to conjugacy
in $\Gamma$ by $g$ and the cycle. For any
conjugacy class $c\in \G$ and each integer $i\geq 1$, 
the number of $i$-cycles in $\sigma$ whose cycle-product lies in $c$
will be denoted by $m_i(c)$. This gives rise to a partition  
$\rho(c)=(1^{m_1 (c)} 2^{m_2 (c)} \ldots )$ for $c \in \G_*$.
Thus we obtain a partition-valued
function $\rho=( \rho (c))_{c \in \G_*} \in {\mathcal P} ( \G_*)$ such that
$\|\rho\|=\sum_{i, c} i m_i(\rho(c)) =n$. This is called the {\em type} of the element
$(g, \sigma)$.  It is known \cite{M2} that two elements
in the same conjugacy class have the same type and there
exists a one-to-one correspondence between the sets
$(\Gn)_*$ and $\mathcal P_n(\G_*)$. 
We will freely say that $\rho$ is the type of the conjugacy class
of $\Gn$. 

Given a class $c$ we denote by $c^{-1}$ the class
$\{x^{-1}| x\in c\}$. 
For each $\rho\in \mathcal P(\G_*)$ we also associate the
partition-valued function 
$$\overline{\rho}=(\rho(c^{-1}))_{c\in \G_*}.$$

Given a partition $\lambda = (1^{m_1} 2^{m_2} \ldots )$,
we denote by 
\[
  z_{\la } = \prod_{i\geq 1}i^{m_i}m_i!
\]
 the order of the centralizer 
of an element of cycle type $\la $ in $S_{|\la |}$.
The order of the centralizer of an element 
$x = (g, \sigma) \in {\Gamma}_n$ of type 
$\rho=( \rho(c))_{ c \in \G_*}$ is given by 
$$
Z_{\rho}=\prod_{c\in \G_*}z_{\rho(c)}\zeta_c^{l(\rho(c))}.
$$

\subsection{Grothendieck ring $R_{\G\times C^{\times}}$}
Let $\Rz$ be the $\mathbb Z$-lattice generated by $\g_i$, $i=0, \dots, r$,
and $R(\G)=\mathbb C\otimes\Rz$ be the space of complex class functions on
 the group
$\G$.
In our previous work on the McKay correspondence and vertex representations
\cite{W, FJW}
we studied the Grothendieck ring $
  \RG = \bigoplus_{n\geq 0} R({\Gamma}_n).$ 
In the quantum case we need to add the ring $R(\mathbb C^{\times})$, 
the space of
characters of $\mathbb C^{\times}=\{t\in\mathbb C|t\neq 0\}$. 

Let $q$ be the 
irreducible character of $\Cs$ that sends $t$ to itself.
Then $R(\mathbb C^{\times})$ is spanned by 
irreducible multiplicative characters $q^n$, $n\in \mathbb Z$, where 
\[
q^n(t)=t^n, \qquad t\in \mathbb C^{\times}.
\]
Thus $R(\mathbb C^{\times})$ is identified with the ring 
$\mathbb C[q, q^{-1}]$, and we have
\[
R(\G\times \mathbb C^{\times})=R(\G)\otimes R(\mathbb C^{\times}).
\]

An elements of $R(\G\times \mathbb C^{\times})$ can be written as a finite sum:
\[
f=\sum_if_i\otimes q^{n_i}, \qquad f_i\in R(\G), n_i\in\mathbb Z.
\]

We can also view $f$ as a function on $\G$ with values in 
the ring of Laurent polynomials $\mathbb C[q, q^{-1}]$. In this case we
will write $f^q$ to indicate the formal variable $q$,
then $f^q(c)=\sum_if_i(c)q^{n_i}\in\mathbb C[q, q^{-1}]$.
As a function on $\G\times \mathbb C^{\times}$, we have
$f(c, t)=\sum_if_i(c)t^{n_i}$.

Denote by $\RGC$ the following direct sum:
\[
\RGC=\bigoplus_{n\geq 0} R(\Gn\times \mathbb C^{\times})
\simeq R_{\G}\otimes \mathbb C[q, q^{-1}].
\]

\subsection{Hopf algebra structure on $\RGC$} 
       The multiplication $m$ in $\mathbb C^{\times}$ and the diagonal map
$\mathbb C^{\times}\stackrel{d}{\longrightarrow}\mathbb C^{\times}\times \mathbb C^{\times}$
induce the Hopf algebra structure on $R(\mathbb C^{\times})$.
\begin{align}\label{E:hopf1}
m_{\mathbb C^{\times}}&: \RC\otimes \RC\stackrel{\cong }{\longrightarrow} 
R(\Cs\times \Cs)
 \stackrel{d^*}{\longrightarrow} R(\Cs),\\ \label{E:hopf2}
\Delta_{\Cs}&: R(\Cs) \stackrel{m^*}{\longrightarrow}
 R(\Cs\times \Cs)
 \stackrel{\cong}{\longrightarrow}
 R(\Cs) \otimes R(\Cs).
\end{align}
In terms of the basis $\{q^n\}$ we have 
\begin{align*}
q^i\cdot q^j&=q^{i+j},\\
\Delta(q^k)&=q^k\otimes q^k,
\end{align*}
where we abbreviate $\Delta_{\Cs}$ by $\Delta$
and follow the convention of writing $a\cdot b=m_{\Cs}(a\otimes b)$.

The antipode $S_{\Cs}$ and the counit $\epsilon_{\Cs}$ are given by
\[
S_{\Cs}(q^n)=q^{-n}, \qquad \epsilon_{\Cs}(q^n)=\delta_{n0}.
\]

We extend the Hopf algebra structures on $R(\Cs)$ and
$\RG$ \cite{Z, M2} into a Hopf algebra structure
on $\RGC$ using a standard procedure in Hopf algebra
\cite{A}.
The multiplication and comultiplication are given by the respective
composition of the following maps:
\begin{align}
 m &: R(\GCn ) \otimes R(\GCm) 
 \stackrel{\cong }{\longrightarrow} R(\GCn \times \GCm)\nonumber\\
 &\qquad \stackrel{1\otimes m_{\Cs}}{\longrightarrow} 
R( {\Gn\times\Gm \times\Cs})\stackrel{Ind\otimes 1}{\longrightarrow}
R( {\Gamma}_{n + m}\times\Cs);\\
\Delta&: R(\GCn ) \stackrel{Res\otimes 1}{\longrightarrow}
 \oplus_{m=0}^nR( {\Gamma}_{n - m} \times \GCm)\nonumber\\
&\qquad \stackrel{1\otimes \Delta_{\Cs}}{\longrightarrow}
\oplus_{m=0}^n R( {\Gamma}_{n - m}\times \GCm\times \Cs)\nonumber\\
&\qquad \stackrel{\cong }{\longrightarrow}
 \oplus_{m=0}^nR( {\Gamma}_{n - m}\times \Cs) \otimes R(\GCm),
\label{E:comult}
\end{align}
where we have used the identification of $R(\Cs\times \Cs)$ with
$R(\Cs)\otimes R(\Cs)$ in (\ref{E:hopf1}-\ref{E:hopf2}). Also
$Ind: R(\Gn \times \Gm ) \longrightarrow R(\G_{n +m})$
denotes the induction functor  
and $Res: R(\Gn ) \longrightarrow
 R( {\G }_{n - m} \times \Gm )$ denotes the restriction functor.

The antipode is given by
\[
S(f(g, t))=f(g^{-1}, t^{-1}), \qquad g\in \G, t\in \Cs.
\]
In particular, $S(\g)(c)=\g(c^{-1})$ for
$\g\in\G^*$. As we mentioned earlier, we may write
$f\in \RGC$ as 
$$f^q(g)=\sum_if_i(g)q^{n_i},$$ 
Then
$S(f^q)(g)=\sum_if_i(g^{-1})q^{-n_i}$.

The counit $\epsilon$ is defined by
\[
\epsilon(R(\GCn ))=0, \qquad\mbox{if \ \ } n\neq 0,
\]
and $\epsilon$ on $R(\Cs)$ is the counit of the Hopf algebra 
$R(\Cs)$.

\section{A weighted bilinear form on $R(\GCn)$}
\label{sect_weight}

\subsection{A standard bilinear form on $\RGC$}
Let $f, g\in R(\GC)$ with $f=\sum_i f_i\otimes q^{n_i}$ and
$g=\sum_i g_i\otimes q^{m_i}$. 
The $\mathbb C[q, q^{-1}]$-valued standard $\mathbb C$-bilinear form
on $R(\GC)$ is defined as
\begin{align*}
\langle f, g \rangle_{\G}^q  
&= \sum_{i,j} \langle f_i, g_j\rangle_{\G}q^{n_i-m_j}\\
&= \sum_{i,j}\sum_{c\in \G_*}\zeta_c^{-1}f_i(c)g_j(c^{-1})q^{n_i-m_j},
\end{align*}
where we recall that
$c^{ -1}$ denotes the conjugacy class 
$\{ x^{ -1}| x \in c \}$ of $\G$, and $\zeta_{c}$ is the order of the 
centralizer of the class $c$ in $\G$. Sometimes we will also view the bilinear form as a function
of $t\in\Cs$:
\begin{equation*}
\langle f, g \rangle_{\G}^q(t)  
= \sum_{c \in \Gamma_*} \zeta_c^{ -1} f(c, t) S(g(c, t)).
\end{equation*}

The following is a direct consequence of the orthogonality
of irreducible characters of $\G$. 
\begin{eqnarray}
  \langle \g_i\otimes q^k, \g_j \otimes q^l\rangle_{\G}^q &= & \delta_{ij}
q^{k-l} , \nonumber \\
  \sum_{ \g \in \G^*} \g (c ') S(\g)( c)
    &= & \delta_{c, c '} \zeta_c, \quad c, c ' \in \G_*.  \label{eq_orth}
\end{eqnarray}

Let $\langle \ \ , \ \ \rangle_{\Gn}^q$ be the 
$\mathbb C[q, q^{-1}]$-valued bilinear form on $R(\GCn)$. 
The $\mathbb C[q, q^{-1}]$-valued standard bilinear form
in $\RGC$ is defined in terms of the bilinear form on $R( \GCn )$
as follows:
\[
\langle u, v \rangle^q
 = \sum_{ n \geq 0} \langle u_n, v_n \rangle_{\Gn}^q,
\]
where
$u = \sum_n u_n$ and $v = \sum_n v_n$ with $u_n, v_n\in R(\Gn\times\Cs)$.

\subsection{A weighted bilinear form on $R(\GC)$}
A class function $\wt \in R(\GC)$ is called {\t self-dual} if for all
$x\in\G, t\in\Cs$
\[
\xi(x, t)=S(\xi(x, t)),
\]
or equivalently $\xi^q(x)=\xi^{q^{-1}}(x^{-1})$.

We fix a self-dual class function $\xi$. 
The tensor product of two representations $\g$ and $\beta$ in $ R(\GC)$
will be denoted by $\g*\beta $.

Let $a_{ij} \in \mathbb C[q, q^{-1}]$ be the (virtual) multiplicity 
of $\g_j$ in $ \wt * \g_i $, i.e.,
\begin{eqnarray}  \label{eq_tens}
 \wt  *  \g_i 
  = \sum_{j =0}^r a_{ij} \g_j.
\end{eqnarray}
We denote by $A^q$ the $ (r +1) \times (r +1)$ matrix 
$ ( a_{ij})_{0 \leq i,j \leq r}$.

Associated to $\xi$ we introduce the following weighted bilinear form
$$
  \langle f, g \rangle_{\wt}^q = \langle \wt * f ,  g \rangle_{\G}^q,
   \quad f, g \in R( \GC).
$$ 
where we use the superscript $q$ to indicate the $q$-dependence.
The superscript $q$ is often omitted if 
the $q$-variable in characters $f$ and $g$ is clear from the context.
The explicit formula of the bilinear form is given as follows.
\begin{eqnarray}
  \langle f, g \rangle_{\wt}^q
   &=&\frac 1{ |\G|} \sum_{x\in \G}\wt^q(x)f^q(x)g^{q^{-1}}(x^{-1})
\nonumber\\ 
& =& \sum_{c \in \G_*} \zeta_c^{ -1} \wt^q(c) f^q(c) g^{q^{-1}}(c^{ -1}),
     \label{eq_twist}
\end{eqnarray}
which is the  average of 
the character $ \wt * f * \overline{g}$ over $\G$.
   
The self-duality of $\wt$ together with (\ref{eq_twist})
implies that 
$$
 a_{ij} = \overline{a_{ji}}, 
$$
i.e. $A^q$ is a hermitian-like matrix with the bar action given by
$\overline{q}=q^{-1}$. 

The orthogonality (\ref{eq_orth}) implies that
\begin{equation}\label{E:qcartan}
a_{ij}=\langle \g_i, \g_j\rangle_{\xi}^q.
\end{equation}

\begin{remark}
 If $\wt$ is the trivial character $\g_0$, then the 
 weighted bilinear form becomes the standard one on $ R( \GC)$.
\end{remark}

\subsection{A weighted bilinear form on $R( \GCn)$}

Let $V$ be a $\G\times \Cs$-module which affords a character $\g$ in $R(\GC)$.
We can decompose $V$ as follows:
\begin{equation*}
V=\bigoplus_iV_i\otimes \mathbb C({k_i}),
\end{equation*}
where $V_i$ is a (virtual) $\G$-module in $R(\G)$ and
$\mathbb C(k_i)$ is the one dimensional $\Cs$-module
afforded by the character $q^{k_i}$. 

The $n$-th outer tensor product 
$V^{ \otimes  n} $ of $V$ can be regarded naturally as 
a representation of the wreath product 
$(\G\times\Cs)_n$ via permutation of the factors
and the usual direct product action. More precisely, note that $\Gn\times \Cs$ can be viewed as a subgroup of $(\G\times \Cs)_n$ by the diagonal inclusion from
$\Cs$ to $(\Cs)^n$:
\[
\Gn\times\Cs\longrightarrow ({\G}^n\times{\Cs}^n)\rtimes S_n=
(\G\times\Cs)_n.
\]
This provides a 
 natural $\Gn\times\Cs$-module structure on $V^{\otimes n}$. We
denote its character by $\eta_n ( \g )$. Explicitly we have
\begin{equation}\label{E:eta-action}
(g, \sigma, t).(v_1\otimes\cdots\otimes v_n)=
(g_{1}, t)v_{\sigma^{-1}(1)}\otimes\cdots 
\otimes (g_{n}, t)v_{\sigma^{-1}(n)},
\end{equation}
where $g=(g_1, \ldots, g_n)\in \G^n$.

Let $\varepsilon_n$ be the (1-dimensional) sign representation 
of $\Gn$ so that $\G^n$ acts trivially while
letting $S_n$ act as a sign representation.
We denote by $\varepsilon_n ( \g ) \in R(\GCn)$
the character of the tensor product
of $\varepsilon_n\otimes 1$ and $V^{\otimes n}$.

The weighted bilinear form on $R( \GCn)$ is now defined by
$$
  \langle  f, g\rangle_{\wt, \Gn }^q =
   \langle \eta_n (\wt ) * f, g \rangle_{\Gn}^q ,
   \quad f, g \in R( \GCn).
$$
We shall see in Corollary~\ref{cor_char} that 
$\eta_n (\wt)$ is self-dual if the class function $\wt$ is invariant under the antipode $S$. In such a case the matrix of the 
bilinear form $\langle \ , \ \rangle_{\wt}^q$ is equal
to its adjoint (transpose and bar action).

We can naturally extend   $\eta_n$ to 
a map from $R(\G)\otimes q^k$ to $R(\GC)$ as in the classical case (cf. \cite{W}).
In particular, if $\beta$ and $\g $ are characters of 
representations $V$ and $W$ of $\G$ respectively, then
\begin{align}  \nonumber
&\eta_n (\beta\otimes q^k +\g\otimes q^l) \\      \label{eq_virt}
  =\sum_{m =0}^n &Ind_{\G_{n -m}\times \Cs \times \Gm\times \Cs }^{\GCn}
   [ \eta_{n -m} (\beta\otimes q^k) \otimes \eta_m (\g\otimes q^l) ],\\
  &\eta_n (\beta\otimes q^k - \g\otimes q^l) \nonumber\\
  =\sum_{m =0}^n &( -1)^m Ind_{\G_{n -m}\times \Cs \times \Gm\times 
\Cs }^{\GCn}
   [ \eta_{n -m} (\beta\otimes q^k) \otimes \varepsilon_m (\g\otimes q^l) ] .  \label{eq_virt'}
\end{align}

On $\RGC = \bigoplus_{n} R(\GCn)$ the weighted 
bilinear form is given by
\[
\langle u, v \rangle_{\wt}^q
 = \sum_{ n \geq 0} \langle u_n, v_n \rangle_{\wt, \Gn }^q
\]
where
$u = \sum_n u_n$ and $v = \sum_n v_n$ with $u_n, v_n\in R(\Gn\times\Cs)$.

 The bilinear form $\langle \ , \ \rangle_{\wt}^q$ on
$\RGC$
is $\mathbb C$-bilinear and takes values in $\mathbb C[q, q^{-1}]$.
When $n =1$, it reduces
 to the weighted bilinear form defined on $R(\GC)$.

We will often omit the superscript $q$ and use the notation $\langle \ , \ \rangle_{\xi}$ for the weighted bilinear form on $R_{\GC}$.

\section{Quantum McKay weights}\label{sect_mckay}

\subsection{Quantum McKay correspondence}
Let $d_i = \g_i (c^0)$ be the dimension of 
the irreducible representation of $\G$
corresponding to the character $\g_i$.

The following generalizes a result of McKay \cite{Mc}.

\begin{proposition} \label{P:qMc}
    For each class $c\in \Gamma_*$ the column vector
 $$
   v(c)= ( \g_0 (c), \g_1 (c), \ldots, \g_r (c) )^t 
   $$
 is an eigenvector of the $(r+1)\times(r+1)$-matrix
 $A^q=(\langle\g_i, \g_j\rangle_{\xi}^q)$
 with eigenvalue $ \wt^q (c)$. 
 In particular $(d_0, d_1, \ldots, d_r)$ is an eigenvector
 of $A^q$ with eigenvalue $\wt^q(c^0)$.
\end{proposition}
\begin{proof} We compute directly that
\begin{align*}
\sum_{k=0}^r\langle\g_i, \g_k\rangle_{\xi}^q\g_k(c)&=\sum_k\sum_{c'\in\Gamma_*}
\zeta_{c'}^{-1}\xi^q(c')\g_i(c')\g_k({c'}^{-1})\g_k(c)\\
&=\sum_{c'\in\Gamma_*}\zeta_{c'}^{-1}\xi^q(c')\g_i(c')\sum_k
\g_k({c'}^{-1})\g_k(c)\\
&=\sum_{c'\in\Gamma_*}\zeta_{c'}^{-1}\xi^q(c')\g_i(c')\zeta_c\delta_{cc'}\\
&=\xi^q(c)\g_i(c).
\end{align*}
\end{proof}

  Let $\pi$ be an irreducible faithful 
representation $\pi$ of $\G$ of
dimension $d$. For each integer $n$ we
define the $q$-integer $[n]$ that can be viewed as a 
character of $\Cs$ by
\begin{equation*}
[n]=\frac{q^n-q^{-n}}{q-q^{-1}}=q^{n-1}+q^{n-3}+\cdots +q^{-n+1}.
\end{equation*}
We take the following
special class function 
\begin{equation}\label{E:weight}
  \wt=\g_0\otimes [d] - \pi\otimes 1_{\Cs},
\end{equation}
where we have also used the symbol $\pi$ for the corresponding
character, and $1_{\Cs}=q^0$ is the trivial character
of $\Cs$.

\begin{proposition} \label{P:nondeg}
The weighted bilinear form associated to 
(\ref{E:weight}) is
non-degenerate. If $\pi$ is an embedding of $\G$ 
into $SU_d$ and
$t\neq 1$ is a nonnegative real number, 
then the weighted bilinear form evaluated on $t$ is positive definite.
\end{proposition}
\begin{proof} A simple fact of finite group theory says that
\[
\langle f, f\rangle_{\pi}\leq d\langle f, f\rangle.
\]
Assume that $t\in\mathbb R_+$. Observe that
\[
t^{d-1}+t^{d-3}+\cdots+t^{-d+1}\geq d
\]
and the equality holds if and only if $t=1$.

Let $A^q=(\langle\g_i, \g_j\rangle_{\wt}^q)=(a_{ij})$. 
Note
that for any faithful representation $\pi$ of $\G$ we have that
\[
\pi*\g_i=\sum_{j}c_{ij}\g_j, \qquad c_{ij}\in\mathbb N.
\]
Then it follows that 
\begin{align*}
\langle \g_i, \g_j \rangle_{\wt}^q(t)&=(t^{d-1}+t^{d-3}+\cdots+t^{-d+1})
\langle \g_i, \g_j \rangle -\langle \g_i, \g_j \rangle_{\pi}\\
&=[d](t)\delta_{ij}-c_{ij}=A^{1}+([d](t)-d)I.
\end{align*}
According to Steinberg (see e.g. \cite{FJW}),
$A^1$ is positive semi-definite which generalizes
McKay's observation in the case of $d=2$.
This implies that the eigenvalues of
$A^q$ are $\geq [d](t)-d\geq 0$. Thus the matrix $A^q(t)$ is positive-definite
when $t>0$ and $t\neq 1$.
\end{proof}

We remark that when $|t|=1$ and $t$ is close to $1$,
the signature of $A^q(t)$ is $(-1, 1, \ldots, 1)$ due to
$[d](t)\leq d$.

\begin{remark} The matrix $A^1$ is integral, and
the entries of $A^q$ are the $q$-numbers of the corresponding
entries in $A^1$ when $r\geq 2$.
\end{remark} 

\subsection{Two quantum McKay weights} \label{S:Mcweights}
Let $\G$ is a finite subgroup of $SU_2$ and we introduce
the first distinguished self-dual class function
$$\xi=\g_0\otimes (q+q^{-1})-\pi\otimes 1_{\Cs},
$$
where $\pi$ is the character of the embedding of $\G$ in $SU_2$. 

The matrix of the
weighted bilinear form $\langle \ , \ \rangle_{\xi}$ (cf. (\ref{E:qcartan}))
has the following entries:
\begin{equation} \label{E:qcartan1}
a_{ij}=\begin{cases} q+q^{-1}, & \mbox{if } i=j,\\
-1, & \mbox{if $\langle \g_i, \g_j\rangle_{\xi}^1=-1$,}\\
-2, & \mbox{if $\langle \g_i, \g_j\rangle_{\xi}^1=-2$ and 
$\G={\mathbb Z}/2{\mathbb Z}$}\\
0, & \mbox{otherwise.}
\end{cases}
\end{equation}

In particular when $q=1$ the matrix $(a_{ij}^1)$ coincides with the
extended Cartan matrix of ADE type according to the five classes of
finite subgroups of $SU_2$: 
the cyclic, binary dihedral, tetrahedral, octahedral, and icosahedral 
groups.
McKay \cite{Mc} gave a direct correspondence between a finite subgroup
of $SU_2$ and the affine Dynkin diagram $D$ of ADE type. 
Each irreducible character $\g_i$ corresponds to  a 
vertex of $D$, and the number of edges between $\g_i$ and $\g_j$ ($i\neq j)$
is
equal to $|\langle \g_i, \g_j\rangle_{\xi}^1|$, where
$\langle \g_i, \g_j\rangle_{\xi}^1=a_{ij}^1$
are the entries of matrix $A^1$ of the weighted bilinear form
$\langle \ , \ \rangle_{\xi}^1$.
For this reason we will call our matrix 
$A^q=(a_{ij})=(\langle\g_i, \g_j\rangle_{\xi}^q)$ 
the quantum Cartan matrix.

Let $\G$ be the cyclic subgroup of $SU_2$ of order $r+1$. 
We can introduce the second deformation
parameter in the quantum Cartan matrix.
Let $\g_i (i=0, \ldots, r)$ be the full set of 
irreducible characters of $\Cs$ such that $\g_i*\g_j
=\g_{i+j\mod r+1}$. The embedding of
$\G$ in $SU_2$ is given by
$\pi=\g_1+\g_r$. 

For $p=q^k\in R(\Cs)$ we let
$$
\xi=\xi^{q, p}
=\g_0\otimes (q+q^{-1})-(\g_1\otimes p+\g_{r}\otimes p^{-1}).
$$
When $p=1$ the second choice reduces to the first choice in type $A$. 
This class function
is self-dual since 
$S(\g_i)=\g_{r+1-i}, i=0, 1, \ldots, r$
and $S(q)=q^{-1}, S(p)=p^{-1}$. 

It is easy to see that  
\begin{equation}\label{E:qcartan2}
a_{ij}(q, p)=\langle\g_i, \g_j\rangle_{\wt}^{q, p}
=[2]\delta_{ij}-p\delta_{i+1, j}-p^{-1}\delta_{i-1, j}.
\end{equation} 
Thus the matrix of the
weighted bilinear form $\langle \ , \ \rangle_{\xi}$ (cf. (\ref{E:qcartan}))
has the following form. 
\begin{equation}\label{E:2varCartan}
\begin{pmatrix}
q+q^{-1} & -p & 0 & \cdots & -p^{-1} \\
-p^{-1} & q+q^{-1} & -p & \cdots & 0 \\
0 & -p^{-1} & q+q^{-1} & \cdots & 0 \\
\cdots & \cdots & \cdots & \cdots & \cdots \\
-p & 0 & \cdots & -p^{-1} & q+q^{-1}
\end{pmatrix}, \qquad \mbox{if } r\geq 2,
\end{equation}
or
\begin{equation}
\begin{pmatrix}
q+q^{-1} & -p-p^{-1}\\
-p-p^{-1} & q+q^{-1}
\end{pmatrix}, \qquad  \mbox{if } r=1.
\end{equation}

Note that 
 when $\G=1$, the matrix of the bilinear form $\langle \ \ , \  \ 
\rangle_{\wt}^{q, p}$ is $q+q^{-1}-p-p^{-1}$, which is degenerate 
when $q=p^{\pm 1}$. 

We will call this matrix the $(q, p)$-Cartan matrix (of type A).
The self-duality of $\wt^{q, p}$ transforms into the
condition that the $(q, p)$-Cartan matrix is $*$-invariant, where the
$*$ action is the composition of transpose and bar action. Namely, 
$a_{ij}(q, p)=a_{ji}(q^{-1}, p^{-1})$.

\begin{proposition} If $p\neq q^{\pm 1}$, then the bilinear form 
$\langle \ , \ \rangle_{\wt}^{q, p}$
is non-degen\-erate. If $p=q^{\pm 1}$, the bilinear form $\langle \ , \ \rangle_{\wt}^{q, p}$ is degenerate of rank $r$.
\end{proposition}
\begin{proof} Let $A^{q, p}=(\langle \g_i, \g_j\rangle_{\xi}^{q, p})$
be the matrix of the bilinear form $\langle \ , \ \rangle_{\wt}^{q, p}$
and let $\omega$ be a $(r+1)$-th root of unity. Then
$\g_i(c^j)=\omega^{ij}$ and $\g_i*\g_j=\g_{i+j}$. From this and 
Proposition \ref{P:qMc} we see that as a matrix over $\mathbb C[q, q^{-1}]$
the eigenvalues of $A^{p, q}$ are
$q+q^{-1}-\omega^i p-\omega^{-i} p^{-1}$, $i=0, \ldots, r$.
The function $q+q^{-1}-\omega^i p-\omega^{-i}p^{-1}\in R(\Cs)$ is non-zero except when $i=0$ and $q=p^{\pm 1}$.
\end{proof}

\section{Quantum Heisenberg algebras and $\Gn$} \label{sect_heis}
\subsection{Heisenberg algebra $\hg $}
Let $\hg $ be the infinite dimensional 
Heisenberg algebra over $\mathbb C[q, q^{-1}]$, associated with 
$\Gamma$ and $\wt\in R(\GC)$,
with generators $a_m(c), c\in\G_*, m\in \mathbb Z$ 
and a central element $C$ subject to the following commutation relations:
\begin{equation}  \label{eq_heis}
[a_m(c^{-1}), a_n(c')]
 = m \delta_{m, -n} \delta_{c, c'}\zeta_{c}\xi_{q^m}(c)C, 
 \quad c, c ' \in \G_*.
\end{equation}

For $m\in \mathbb Z, \g\in \G^*$ and $k\in \mathbb Z$ we define
\begin{equation*}
  a_m( \g\otimes q^k ) 
    = \sum_{c \in \G_*} \zeta_c^{ -1}
  {\g } (c) a_m(c)q^{mk}                   
\end{equation*}
and then extend it to $R(\GC)$ linearly over $\mathbb C$. Thus
we have for $\g \in R(\GC)$
\begin{equation}\label{eq_real}
  a_m(\g) 
    = \sum_{c \in \G_*} \zeta_c^{ -1}
  {\g_{q^m} } (c) a_m(c).                   
\end{equation}
In particular we have $a_m(\g\otimes q^k)=a_m(\g)q^{mk}$.

It follows immediately from the orthogonality (\ref{eq_orth}) of the irreducible
characters of $\G $ that for each $c\in \G_*$ 
$$
 a_{ m}( c) = \sum_{ \g\in \G^*} S(\gamma(c)) a_m( \g ).
$$
Note that this formula is also valid 
if the summation runs through $\G^*\otimes q^k$ with a fixed $k$.

\begin{proposition}  \label{prop_orth} The Heisenberg
 algebra $\hg$ has a new basis given by
$a_{n}(\g )$ and $C$ ($ n\in \mathbb Z, \g \in \G^*$) 
over $\mathbb C[q, q^{-1}]$ 
 with the following relations:
 \begin{equation} \label{E:heisen}
  [ a_m(\g), {a_n( \g' )}]=m\delta_{m, -n}\langle\g, \g'
\rangle_{\xi}^{q^m}C.
 \end{equation}
\end{proposition}

\begin{proof} This is proved by a direct computation using
Eqns. (\ref{eq_heis}),
 (\ref{eq_twist}) and (\ref{eq_orth}).
 \begin{eqnarray*}
  [ a_m(\g), a_n( \g')]
  & =& \sum_{c, c'\in \G_*} \zeta_c^{-1}\zeta_{c'}^{-1}
{ \g  } (c)\g'(c') 
       [ a_m( c), a_n ({c'})]          \\
& =& m \delta_{m, -n}
         \sum_{c, c'\in \G_*}\zeta_{c}^{-1}\zeta_{c'}^{-1}
{ \g  } (c)\g'(c') \delta_{c^{-1}, c'}\zeta_c\xi_{q^m}(c)C    \\
  & =& m \delta_{m, -n}
         \sum_{c\in \G_*}\zeta_{c}^{-1}
{ \g  } (c)\g'(c^{-1}) \xi_{q^m}(c)C   \\
  & =& m \delta_{m, -n} \langle\g, \g'\rangle_{\xi}^{q^m}C. 
\end{eqnarray*}
\end{proof}

\subsection{Action of $\hg$ on the Space $\SGC$}
 Let $\SGC $ be the symmetric algebra 
generated by $a_{-n}(\g), n \in \mathbb N,
\g\in \Gamma_*$ over $\mathbb C[q, q^{-1}]$.
We define $a_{-n}(\g\otimes q^k)=a_{-n}(\g)q^{-kn}$
and the natural degree operator on the space $\SGC $ by
$$
  \deg (a_{ -n}( \g\otimes q^k)) = n 
$$ 
which makes $\SGC $ into a $\mathbb Z_+$-graded algebra. 

The space $\SGC$ affords a natural realization of the
Heisenberg algebra $\hg$ with $C=1$. 
Since $a_{-n}(\g\otimes q^k)=q^{-nk}a_{-n}(\g)$, it is enough to describe
the action for $a_{-n}(\g)$.
The central element $C$ acts as the identity operator.
For $n>0$, $a_{-n}( \g)$ act as     
multiplication operators on $ \SGC $.
The element $a_n (\g), n \geq 0$ acts as a differential operator through
contraction:
\begin{eqnarray*}
 & a_n (\g). a_{-n_1}( \alpha_1) a_{-n_2} (\alpha_2)
    \ldots a_{-n_k}( \alpha_k)   \\
 &= \sum_{i =1}^k 
 \langle \g , \alpha_i \rangle_{\wt}^{q^n}
   a_{-n_1}( \alpha_1) a_{-n_2}(\alpha_2) \ldots 
  \check{a}_{-n_i}( \alpha_i) \ldots a_{-n_k}(\alpha_k )  .
\end{eqnarray*}
Here $n_i > 0, \alpha_i \in R(\G)$ for $i =1, \ldots , k$,
and $\check{a}_{-n_i}( \alpha_i)$ means the very term 
is deleted. 
In this case  $\SGC$ is an irreducible
representation of $\hg$ with the unit $1$ as the highest weight vector.

\subsection{The bilinear form on $\SGC $}
       As a $\hg$-module, the space $ \SGC $ admits a bilinear form
$\langle \ ,  \ \rangle_{\wt } '$ over $\mathbb C[q, q^{-1}]$ characterized by
\begin{align} \nonumber
\langle 1, 1\rangle_{\xi}'&=1,\\   \label{eq_form}
\langle au, v\rangle_{\xi}'&=\langle u, a^*v\rangle_{\xi}', \qquad a\in \hg,
\end{align}
with the adjoint map $*$ on $\hg$ given by
\begin{equation}  \label{eq_hermit}
a_n(\g\otimes q^k)^* = a_{-n}(\g\otimes q^{k}), 
\qquad n\in \mathbb Z. 
\end{equation}
Note that the adjoint map $*$ is a $\mathbb C$-linear 
anti-homomorphism of $\hg$, and $q^*=\overline q$.
We still use the same symbol $*$ to denote the
hermitian-like dual, since it clearly generalizes the $*$-action
on the deformed Cartan matrix (\ref{E:2varCartan}).

  For any partition $\la =( \la_1, \la_2, \dots)$ and $\g \in \G^*$, 
we define 
$$
  a_{-\la}( \g) = a_{-\la_1}( \g)a_{ - \la_2}( \g) \dots .
$$
For $\rho = ( \rho (\g) )_{ \g \in \G^*} \in {\mathcal P}(\G^* )$, 
we define 
$$
  a_{ - \rho\otimes q^k } = q^{-k\|\rho\|}\prod_{\g \in \G^*}  a_{ - \rho (\g)}(\g).
$$
It is clear that for a fixed $k\in \mathbb Z$ the
elements $a_{ - \rho\otimes q^k},
\rho \in {\mathcal P}(\G^* )$ form a basis of 
$\SGC $ over $\mathbb C[q, q^{-1}]$. 

Given a partition $ \la = ( \la_1, \la_2, \ldots )$
and $c \in \G_*$, we define
\begin{eqnarray*}
  a_{ - \la } (c\otimes q^k ) & =& q^{-k|\la|}
a_{ - \la_1}(c) a_{ - \la_2 } (c) \ldots, \\
\end{eqnarray*}
For any $\rho = ( \rho (c) )_{ c \in \G_* } \in 
 \mathcal P ( \G_* )$ and $k\in \mathbb Z$, we define
\begin{equation*}
  a_{- \rho\otimes q^k}' = q^{-k\|\rho\|}
\prod_{ c \in \G_*} a_{ - \rho (c)} (c).
\end{equation*}

It follows from
Proposition~\ref{prop_orth} that
\begin{eqnarray}  \label{eq_inner}
  \langle a_{ - \rho\otimes q^k}', {a_{ - \overline{\sigma}\otimes q^{l}}'} 
\rangle_{\wt }'
  = \delta_{\rho, \sigma}q^{\|\rho\|(l-k)}
    Z_{\rho} \prod_{c \in \G_*} \prod_{i\geq 1}\wt_{q^i}(c)^{m_i(\rho (c))},
 \end{eqnarray}
where $\rho, \sigma \in \mathcal P (\G_*)$. Note that 
$S(a_{-\rho\otimes q^k}')=a_{-\overline{\rho}\otimes q^{-k}}'$, where we recall that
$\overline{\rho}
\in {\mathcal P}(\G_* )$
is the partition-valued function given by $c\mapsto \rho(c^{-1})$, $c\in \G$.

\section{The characteristic map as an isometry} 
\label{sect_isom}
\subsection{The characteristic map $\ch$}
  Let $\Psi : \Gn \rightarrow \SGC$ be the map defined
by $\Psi (x) = a_{ - \rho}'$ if $x \in \Gn$ is of type $\rho$.

      We define a $\mathbb C$-linear map 
$ch: \RGC \longrightarrow \SGC$ by letting
\begin{align} \nonumber
 ch (f ) 
 &= \langle f, \Psi \rangle_{\Gn}\\
 &= \sum_{\rho \in \mathcal P(\G_*)} Z_{\rho}^{-1} S(f(\rho)) 
 a_{ - \rho}',  \label{E:ch}
\end{align}
where $f(\rho)\in \mathbb C[q, q^{-1}]$ is the 
value of $f$  at the elements of type $\rho$.
The map $ch $ is called the {\em characteristic map}. This generalizes
the definition of the characteristic map in the classical setting
(cf. \cite{M2, FJW}).

The space $\SGC$ can also be interpreted as follows.
The element $a_{ -n } (\g ), n >0 , \g \in \G^*$ is identified as
the $n$-th power sum in a sequence of variables 
$ y_{ g } = ( y_{i\g } )_{i \geq 1}$. By the commutativity among
$a_{-n}(\g)$ ($\g\in\G^*, n>0$) and dimension counting it is clear that 
the space $\SGC$ is isomorphic with the space $\LG$ of 
symmetric functions indexed by $ \G^*$ tensored with
$\mathbb C[q, q^{-1}]$ (cf. \cite{M2}).

Denote by $c_n (c \in \G_*)$ the conjugacy class in $\Gn$
of elements $(x, s) \in \Gn$ such that $s$ is an $n$-cycle
and $x \in c$. 
Denote by $\sigma_n (c\otimes q^k )$ the class function on $\Gn\times\Cs$ which takes 
values $n \zeta_ct^{-nk}$ (i.e. the order of the centralizer of
an element in the class $c_n$ times $t^{-nk}$) 
on elements in the class $c_n\times t$
and $0$ elsewhere. 
For
$\rho = \{ m_r (c) \}_{r \geq 1, c \in \G_*} 
\in \mathcal P_n (\G^*)$ and $k\in\mathbb Z$, 
$$\sigma_{\rho\otimes q^k} = q^{-nk}
\prod_{r \geq 1, c \in \G_*} \sigma_r (c)^{m_r (c)}
$$
is the class function on $\Gn\times\Cs$ which takes value
$Z_{\rho}t^{-nk}$ on the conjugacy class of type $\rho\times t$ and
$0$ elsewhere. Given $\g\in \G^*$ and $k\in \mathbb Z$, we
denote by $\sigma_n (\g \otimes q^k)$ the class function on $\GCn$ which takes 
values $n \g (c)t^{-nk}$ on elements in the class 
$c_n\times t (c \in \G_*)$ 
and $0$ elsewhere. 

\begin{lemma}  \label{lem_isom}
  The map $ch$ sends $\sigma_{\rho\otimes q^k}$ to $a_{ - \rho\otimes q^k} '$.
 In particular, it sends  $\sigma_n(\g\otimes q^k )$ to $a_{ -n} ( \g\otimes q^k )$ in $\SGC$.
\end{lemma}
\begin{proof} This is verified by the definition
of $ch$ (\ref{E:ch})
and the character values of $\sigma_n$ defined above.
\end{proof}

\begin{proposition} 
Given $\g\in \G^*$, the
 character value of 
$\eta_n(\g\otimes q^k)$ on the conjugacy class $c_{\rho}$
of type $\rho=(\rho(c))_{c\in\G_*}$ 
is given by
\begin{equation} \label{E:charval}
\eta_n(\g\otimes q^k)(c_{\rho})
=\prod_{c\in \G_*}\g(c)^{l(\rho(c))}q^{nk}.
\end{equation}
In particular, we have $\eta_n(\g\otimes q^k)=\eta_n(\g)q^{nk}$.
\end{proposition}
\begin{proof}
We first let $(g, \sigma)$ be an element of $\Gn$ such that $\sigma$ is a cycle of length $n$, say $\sigma=(12\cdots n)$. Let $\{e_i\}$ be a basis of $V$, and $\g\otimes q^k$ is afforded
by the action: $(h, t)e_j=\sum_{i}c_{ij}(h)t^ke_i$, where
$h\in \G$. We then have 
\begin{align*}
&(g, \sigma, t).(e_{j_1}\otimes e_{j_2}\otimes \cdots\otimes e_{j_n})\\
&=(g_1, t)e_{j_n}\otimes (g_2, t)e_{j_1} \otimes \cdots\otimes (g_n, t)e_{j_{n-1}}\\
&=\sum_{i_1,\ldots, i_n}t^{kn} c_{i_{n}j_n}(g_1)c_{i_1j_1}(g_2)\cdots
c_{i_{n-1}j_{n-1}}(g_n)
e_{i_n}\otimes e_{i_1}\cdots\otimes e_{i_{n-1}}.
\end{align*}

It follows that
\begin{eqnarray*}
  \eta_n (\g\otimes q^k) (c_{\rho}, t) &=& \mbox{trace }(g, \sigma, t)\\
&=& \sum_{j_1, \ldots, j_n} t^{kn}c_{j_1j_n}(g_1)c_{j_2j_1}(g_2)\cdots
c_{j_nj_{n-1}}(g_n)\\
& =& \mbox{trace } t^{kn}a(g_n) a(g_{n-1}) \ldots a(g_1)   \\
  & =& \mbox{trace } t^{kn}a(g_n g_{n -1} \ldots g_1) = \g (c)q^{kn}(t).
 \end{eqnarray*}

 Given $x\times y \in \Gn$ where
 $x \in \G_r$ and $y \in \G_{n -r}$, by (\ref{E:eta-action}) we clearly have 
$$
\eta_n (\g\otimes q^k) (x\times y, t) = \eta_n (\g\otimes q^k) 
(x, t) 
\eta_n (\g\otimes q^k)(y, t).
$$
 This immediately implies the formula.
\end{proof}

A similar argument gives that
\begin{equation}
  \varepsilon_n (\g\otimes q^k ) ( x, t)
    = (-1)^n  \prod_{c\in \G_*} 
( - \g (c))^{l(\rho(c))}t^{nk},
   \label{eq_signterm}
 \end{equation}
where $x$ is any element in the conjugacy class of type 
$\rho=(\rho(c))_{c\in\G^*}$.

Formula (\ref{E:charval}) is equivalent to the following:
\begin{equation} \label{E:charval''}
\eta_n(\g\otimes q^k)(c_{\rho}, t)=\prod_{c\in \G_*}
\prod_{i\geq 1}
(\g\otimes q^k)(c, t^{i})^{m_i(\rho(c))}.
\end{equation}

The following result allows us to extend the map from
$\g\in\G^*$ to $R(\G_n)$.

\begin{proposition}  \label{prop_exp}
   For any $\g \in R(\G)$, we have
\begin{align}
 \sum\limits_{n \ge 0}  \ch ( \eta_n( \g\otimes q^k ) ) z^n 
  &= \exp \Biggl( \sum_{ n \ge 1} 
      \frac 1n \, a_{-n}(\g)(q^{-k}z)^n \Biggr), \label{eq_exp} \\
 \sum\limits_{n \ge 0}  \ch ( \varepsilon_n( \g\otimes q^k )  ) z^n
  &= \exp \Biggl( \sum_{ n \ge 1} 
      ( -1)^{ n -1} \frac 1n \, a_{-n}(\g)(q^{-k}z)^n \Biggr). 
\label{eq_sign}
\end{align}
\end{proposition}

\begin{proof} It follows from definition of ch (\ref{E:ch}) and 
(\ref{E:charval''}) that
 
\begin{eqnarray*}
 &&\sum\limits_{n \ge 0}  \ch ( \eta_n( \g\otimes q^k ) ) z^n\\ 
  &= & \sum_{\rho} Z_{\rho}^{ -1} 
         \prod_{c\in \G_*}\prod_{i\geq 1}S(\g_{q^{ik}} (c)^{m_i(\rho(c))})
          a_{ -\rho (c) } z^{|| \rho||}q^{-||\rho||}                  \\
 &= & \sum_{\rho} Z_{\rho}^{ -1} 
         \prod_{c\in \G_*}\g (c)^{l(\rho(c))}
          a_{ -\rho (c) } (q^{-k}z)^{|| \rho||}                  \\
  &= & \prod_{c\in \G_*} \Bigl ( \sum_{\lambda }
         (\zeta_c^{ -1}\g (c) )^{l (\lambda)}
         z_{\lambda}^{-1} a_{- \lambda} (c) (q^{-k}z)^{|\lambda|} \Bigr )    \\
  &= & \exp \Biggl  ( \sum\limits_{ n \geq 1}
         \frac1n \sum\limits_{c \in \G_*}
           \zeta_c^{ -1} \g(c) a_{-n} (c) (q^{-k}z)^n \Biggl )      \\
  &= & \exp \Biggl( \sum_{ n \ge 1} 
         \frac 1n \, a_{-n}(\g )(q^{-k}z)^n \Biggr).
 \end{eqnarray*}
 Similarly we can prove (\ref{eq_sign}) using
 the following identity
\begin{eqnarray}
  \varepsilon_n (\g\otimes q^k ) ( x)
    &= &(-1)^n  \prod_{c\in \G_*} \prod_{i\geq 1}
( - \g_{q^{ik}} (c))^{m_i(\rho(c))} \nonumber\\
&= &(-q^k)^n  \prod_{c\in \G_*} \prod_{i\geq 1}
( - \g (c))^{m_i(\rho(c))} \nonumber\\
&=&\varepsilon_n (\g) ( x)q^{nk}.
   \label{eq_signterm'} \nonumber
 \end{eqnarray}

The same argument as in the classical case 
(cf. \cite{FJW}) by using (\ref{eq_virt}) and (\ref{eq_virt'})
will show that
the proposition holds for linear combination of simple
characters such as $\g\otimes q^k-\beta\otimes q^k$, and thus it is true for 
any element $\g\otimes q^k$, where $\g\in  R(\G)$. 
\end{proof}

Comparing components we obtain 
\begin{eqnarray*}
  \ch (\eta_n (\g\otimes q^k ) )
  &=& \sum\limits_{\la}\frac {q^{-nk}}{z_\la }\,
             a_{-\lambda}(\g ), \\
  \ch (\varepsilon_n (\g\otimes q^k ))
  &=& \sum\limits_{\la}\frac {q^{-nk}}{z_\la }\,
             ( -1)^{ | \la | - l ( \la )} a_{-\lambda}(\g ),
\end{eqnarray*}
where the sum runs over all partitions $\lambda$ of $n$. 

\begin{corollary}  \label{cor_char}
   The formula (\ref{E:charval''}) remains valid when
$\g\otimes q^k$ is replaced by any element $\wt \in R(\GC)$.
 In particular $\eta_n (\wt)$ is self-dual provided that $\wt$ is
invariant under the antipode $S$.
\end{corollary}

\subsection{Isometry between $\RGC$ and $\SGC$}
    The symmetric algebra
$\SGC=\SG\otimes \mathbb C[q, q^{-1}]$ 
has the following Hopf algebra structure over $\mathbb C$.
The multiplication is the usual one, and the comultiplication is given by
\begin{align*}
\Delta(q^k)&=q^k\otimes q^k\\
  \Delta ( a_n (\g\otimes q^k ))
&=a_n(\g\otimes q^k)\otimes q^{nk}+q^{nk}\otimes a_n(\g\otimes q^k ),
\end{align*}
where $\g\in\G^*$. The last formula is equivalent to the following:
\begin{equation}\label{E:coprod}
\Delta ( a_n (c\otimes q^k ))
=a_n(c\otimes q^k)\otimes q^{nk}+q^{nk}\otimes a_n
(c\otimes q^k ),
\end{equation}
where $c\in\G_*$.
The antipode is given by
\begin{align*}
S(q^k)&=q^{-k},\\
S(a_n (\g\otimes q^k ))&=-a_n(\g\otimes q^{-k})
\end{align*}
The antipode commutes with the adjoint (dual) map $*$:
\begin{equation}
*^2=S^2=Id, \qquad S*=*S.
\end{equation}

Recall that we have defined a Hopf algebra structure
on $\RGC$ in Sect.~\ref{sect_wreath}. 

\begin{proposition} \label{P:isometry1}
  The characteristic map $ \ch: \RGC \longrightarrow \SGC$
 is an isomorphism of Hopf algebras.
\end{proposition}

\begin{proof} 
It follows immediately from the definition of the comultiplication in 
the both Hopf algebras (cf. (\ref{E:comult}) and (\ref{E:coprod})). 
\end{proof}

\begin{remark} The comultiplication (\ref{E:coprod}) 
is in fact induced from that of the classical case in \cite{FJW}
and only works for $C=1$. 
\end{remark}

\begin{remark} There is another coproduct called Drinfeld comultiplication
$\Delta_D$ on the algebra $\SGC$ adjoined by a central element $q^c$. 
The formula on $\SGC$ at level $c$ is as follows \cite{J2}:
\begin{equation}
\Delta_D(a_n(\g))=a_n(\g)\otimes q^{|n|c/2}+q^{-|n|c/2}\otimes a_n(\g).
\end{equation}
We do not know a conceptual interpretation of the Drinfeld comultiplication
in $\RGC$.
\end{remark}

    Recall that we have defined a bilinear form 
$\langle \  ,  \, \rangle_{\wt }$ on $\RGC$ and
a bilinear form on $\SGC$ denoted by
$\langle \ , \, \rangle_{\wt }'$, where $\wt$ is a self-dual
class function. The following
lemma is immediate from our definition of
$\langle \ , \, \rangle_{\wt }'$ and the
comultiplication $\Delta$.

\begin{lemma}
  The bilinear form $\langle \  ,  \, \rangle_{\wt } '$ on $\SGC$
 can be characterized by the following two properties:
 
 1). $\langle a_{ -n} (\beta\otimes q^k ), a_{ -m} (\g\otimes q^l ) 
\rangle_{\wt}^{'}
  = \delta_{n, m} q^{n(l-k)}\langle \beta , \g  \rangle_{\wt}^{'} ,
  \quad \beta, \g \in \G^*,$ $ k, l\in\mathbb Z.$

 2). $ \langle f g , h \rangle_{\wt}^{'}
       = \langle f \otimes g, \Delta h 
\rangle_{\wt}^{'} ,$
 where $f, g, h \in \SGC $, and the bilinear form on 
$\SGC \otimes \SGC$, is induced
 from $\langle \ , \ \rangle_{\wt}^{'}$ on $\SGC$. 
\end{lemma}

\begin{theorem}  \label{th_isometry}
  The characteristic map is an isometry from the space
 $ (\RGC, \langle \ \ , \ \  \rangle_{\wt })$ to the space
 $ (\SGC, \langle \ \ , \ \  \rangle_{\wt }' )$.
\end{theorem}

\begin{proof}
  By Corollary~\ref{cor_char}, the character
 value of $\eta_n (\wt )$ at an element $x$ of type $\rho$ is
 $$
  \eta_n (\wt ) ( x)= \prod_{c\in \G_*}\prod_{i\geq 1} \wt_{q^i} (c)^{m_i(\rho(c))}.
 $$
 Thus it follows from definition that
 \begin{eqnarray*}
   \langle \sigma_{ \rho\otimes q^k}, \sigma_{  \rho '\otimes q^l} 
\rangle_{\wt }
  & =& \sum_{\mu \in \mathcal P_n (\G_*)}
       Z_{\mu}^{ -1} q^{n(l-k)}\wt_{q} (c_{\mu }) \sigma_{\rho} (c_{\mu})
       \sigma_{\rho '} (c_{\mu})                      \\
  & =& \delta_{\rho, \rho '}
       Z_{\rho}^{ -1}q^{n(l-k)}\wt (c_{\rho}) Z_{\rho}Z_{\rho}  \\
  & =& \delta_{\rho, \rho '}
       Z_{\rho} q^{n(l-k)}
\prod_{c \in \G_*}\prod_{i\geq 1} \wt_{q^i} (c)^{m_i(\rho(c))}.
 \end{eqnarray*}

  By Lemma~\ref{lem_isom} and the formula (\ref{eq_inner}), we 
 see that
 \[
  \langle \sigma_{ \rho\otimes q^k}, \sigma_{  \rho '\otimes q^l} 
\rangle_{\wt }
  = \langle a_{- \rho\otimes q^k}, a_{ - \rho '\otimes q^l} \rangle_{\wt } '
  =  \langle \ch (\sigma_{ \rho\otimes q^k}), 
             \ch (\sigma_{  \rho '\otimes q^l} ) \rangle_{\wt } '.
 \]
 Since $\sigma_{\rho\otimes q^k}, \rho \in \mathcal P(\G_*)$ form
 a $\mathbb C$-basis of $\RGC$, we have shown that 
 $\ch: \RGC \longrightarrow \SGC$ is an isometry.
\end{proof}


>From now on we will not distinguish the bilinear form
$\langle \ , \ \rangle_{\wt}$ on $\RGC$
 from the bilinear form $\langle \ , \ \rangle_{\wt}^{'}$ on $\SGC$.

\section{Quantum vertex operators and $\RGC$}
\label{sect_vertex}
\subsection{Vertex Operators and Heisenberg algebras in
$\FGC$}
 Let $Q$ be an integral lattice with basis $\a_i$, $i=0, 1, \ldots, r$
endowed with a symmetric bilinear form. 
As in the case of $q=1$ (cf. \cite{FK}), we fix a 
$2$-cocycle $\epsilon: Q\times Q \longrightarrow \mathbb C^{\times}$
such that
$$
\ep (\a, \beta)=\ep(\beta, \a)(-1)^{\langle \alpha, \beta \rangle +
        \langle \alpha, \alpha \rangle \langle \beta , \beta \rangle}. 
$$
We remark that the cocycle can be constructed directly by prescribing the 
values of $(\a_i, \a_j)\in \{\pm 1\}$ $(i<j)$.

Let $\xi$ be a self-dual virtual character in
$\RGC$. Recall that the lattice $\Rz$ is a 
$\mathbb Z$-lattice
under  the bilinear form
$\langle \ , \ \rangle_{\wt}^1$, here 
the superscript means $q=1$. For our purpose
we will always associate a $2$-cocycle $\ep$ 
as in the previous subsection to the integral lattice
$(\Rz, \langle \ , \ \rangle_{\wt}^1 )$ (and its sublattices).
 
Let $\mathbb C[\Rz]$ be the group algebra
generated by $e^{\g }$, $\g \in \Rz$.
We introduce two special operators acting on $\mathbb C[ \Rz ]$:
A ($\ep$-twisted) multiplication operator $e^{\alpha}$ defined by
 $$
  e^{\alpha }.e^{\beta } = \ep(\alpha, \beta) e^{\alpha +\beta},
 \quad  \alpha, \beta  \in \Rz,
 $$
and a differentiation operator 
${\partial_{\alpha }}$ given by
\begin{eqnarray*}
 {\partial_{\a}} e^{ \beta} = 
 \langle \a, \beta \rangle_{\wt}^1  e^{ \beta},
 \quad  \alpha, \beta  \in \Rz.
\end{eqnarray*}
These two operators are then extended linearly to the space
\begin{equation}\label{E:fgc}
\FGC = \RGC \otimes \mathbb C[\Rz]
\end{equation} 
by letting them act on the $\RGC$ part trivially.

We define the Hopf algebra structure on $\mathbb C[\Rz]$ 
and extend the Hopf algebra structure from $\RGC$ to $\FGC$as follows.
\begin{equation*}
\Delta(e^{\alpha})=e^{\alpha}\otimes e^{\alpha}, 
\qquad S(e^{\alpha})=e^{-\alpha}.
\end{equation*}

The bilinear form $\langle\ , \ \rangle_{\wt}^q$ on
$\RGC$ is extended
to $\FGC$ by
\begin{equation*}
\langle e^{\alpha}, e^{\beta}\rangle_{\wt}=\delta_{\a,\beta}.
\end{equation*}

With respect to 
this extended bilinear form we have the $*$-action (adjoint action)
 on the 
operators $e^{\alpha}$ and ${\partial}_{\alpha}$:
\begin{equation}
(e^{\alpha})^*=e^{-\alpha}, \qquad 
(z^{{\partial}_{\alpha}})^*=z^{-{\partial}_{\alpha}}.
\end{equation}

For each $k\in\mathbb Z$, we introduce the group theoretic operators \linebreak
$ H_{ \pm n}( \g\otimes q^k ),  
E_{ \pm n} ( \g\otimes q^k), \g \in R(\G), n > 0 $ 
as the following compositions of maps:
\begin{eqnarray*}
  H_{ -n} ( \g\otimes q^k ) &:&
    R ( \GCm ) 
  \stackrel{ \eta_n (\g\otimes q^k) \otimes}{\longrightarrow}
    R ( \GCn ) \otimes R ( \GCm )\\&&
  \stackrel{ {Ind}\otimes m_{\Cs} }{\longrightarrow}
    R ( \G_{n +m}\times\Cs )   \\
  E_{ -n} ( \g\otimes q^k ) &:&
    R ( \GCm ) 
  \stackrel{ \varepsilon_n (\g\otimes q^k) \otimes}{\longrightarrow}
    R ( \GCn ) \otimes R ( \GCm )\\&&
  \stackrel{ {Ind}\otimes m_{\Cs} }{\longrightarrow}
    R ( \G_{n +m}\times\Cs )   \\
  E_n ( \g\otimes q^k ) &:&
    R ( \GCm ) 
   \stackrel{ {Res} }{\longrightarrow}
    R( \Gn) \otimes R( \G_{m -n}\times\Cs)\\&&
   \stackrel{ \langle \varepsilon_n (\g\otimes q^k),
         \cdot \rangle_{\wt } }{\longrightarrow}
    R ( \G_{m -n}\times\Cs) \\
  H_n(\g\otimes q^k ) &:&
    R ( \GCm ) 
   \stackrel{ {Res} }{\longrightarrow}
    R( \Gn) \otimes R( \G_{m -n}\times\Cs)\\&&
   \stackrel{ \langle \eta_n (\g\otimes q^k), 
\cdot \rangle_{\wt} }{\longrightarrow}
    R ( \G_{m -n}\times\Cs) ,
\end{eqnarray*}
where $Res$ and $Ind$ are the restriction and induction functors
in $R_{\G}=\bigoplus_{n\geq 0}R(\G_n)$.

\medskip

We introduce their generating functions in a formal variable $z$:
\begin{eqnarray*}
  H_{\pm} (\g\otimes q^k, z) &=& \sum_{ n\geq 0} H_{ \mp n} 
( \g\otimes q^k ) z^{\pm n}, \\
  E_{\pm} (\g\otimes q^k, z) &=& \sum_{ n\geq 0} E_{ \mp n} 
( \g\otimes q^k )( -z)^{\pm n}.
\end{eqnarray*}
  We now define 
the vertex operators $Y_n ^{\pm}(\g\otimes q^l, k)$ , 
$\g\in \G^*$, $k, l\in \mathbb Z$,
$n \in {\mathbb Z} + \langle \g, \g \rangle_{ \wt}^1 /2$ 
as follows.     
\begin{align}  \nonumber
 Y ^{+}( \g\otimes q^l, k, z)
  & = \sum\limits_{n \in 
  {\mathbb Z} + \langle \g, \g \rangle_{ \wt}^1 /2}
 Y_n^{+}( \gamma\otimes q^l, k)
z^{ -n - \langle \g, \g \rangle_{ \wt}^1 /2}       
\\
  & =  H_+ (\g\otimes q^l, z) E_- (\g\otimes q^{l-k} , z) e^{ \g} 
(q^{-l}z)^{ \partial_{ \g}}, \label{eq_vo}      
\end{align}
\begin{align}\nonumber
Y^-(\g\otimes q^l, k, z)&=(Y^+(\g\otimes q^l, k, z^{-1}))^*\\ \nonumber
&=\sum\limits_{n \in 
  {\mathbb Z} + \langle \g, \g \rangle_{ \wt}^1 /2}
 Y_n^{-}( \gamma\otimes q^l, k)
z^{ -n - \langle \g, \g \rangle_{ \wt}^1 /2}\\ 
&=E_+ (\g\otimes q^{l-k}, z) 
H_- (\g\otimes q^l , z) e^{ -\g} 
(q^{-l}z)^{ -\partial_{ \g}}.
\end{align}

One easily sees that the operators $Y_n^{\pm} (\g\otimes q^l, k)$
are well-defined operators acting on the space
$\FGC$.
 
We extend the $\mathbb Z_+$-gradation 
on $\RGC$ to a $\frac12\langle \g , \g \rangle_{\wt}^1 +
 \mathbb Z_+$-gradation on
$\FGC$ by letting
\begin{eqnarray*}
  \deg  a_{ -n} (\g\otimes q^k ) = n , \quad
  \deg e^{\g } = \frac12 \langle \g , \g \rangle_{\wt}^1 .
\end{eqnarray*}

We denote by $\RGGC$ the subalgebra of $\RGC$
excluding the generators
$a_n(\g_0)$, $n\in \mathbb Z^{\times}$. 
The bilinear form $\langle \ , \ \rangle_{\xi}$ on 
$$\FGGC=\RGGC\otimes\Rzz
$$ 
will be the restriction of $\langle \ , \ \rangle_{\xi}$
on $\FGC$ to $\FGGC$.
In the case of the second choice of $\xi$ and $p=q^{\pm 1}$, the 
Fock space $\FGGC$ can also be obtained as
the quotient of $\FGC$ modulo the radical
of $\langle \ , \ \rangle_{\xi}$.

We define $ \widetilde{a}_{ -n} (\gamma\otimes q^k), n >0$ to be a map 
from $\RGC$ to itself by the following composition
\begin{align*}
  R (\GCm) \stackrel{ \sigma_n ( \g\otimes q^k ) \otimes }{\longrightarrow}&
  R(\GCn) \otimes R (\GCm) \\
& \stackrel{{Ind\otimes m_{\Cs}} }{\longrightarrow}
  R ( {\Gamma}_{n +m}\times \Cs).
\end{align*}
We also define $ \widetilde{a}_{ n} (\gamma\otimes q^k), n >0$ to be 
a map from $\RGC$ 
to itself
as the composition
\begin{align*}
  R (\GCm)  \stackrel{ Res \otimes 1}{\longrightarrow}&
   R(\GCn)\otimes R ( {\G }_{m -n}\times \Cs)\\
 &\stackrel{ \langle \sigma_n ( \g\otimes q^k), \cdot \rangle_{\wt}^q}{\longrightarrow}
 R ( {\G }_{m -n}\times\Cs).
\end{align*}

\begin{proposition} The operators $\widetilde{a}_{n} (\gamma)$, 
$\g\in \G^*, n\in\mathbb Z^{\times}$ satisfy the
Heisenberg algebra relations (\ref{eq_heis}) with $C=1$.
\end{proposition}
\begin{proof} This is similarly proved as for the classical
setting in \cite{W}.
\end{proof}

\subsection{Group theoretic interpretation of vertex operators}
To compare the vertex operators $Y^{\pm}(\g\otimes q^l, k, z)$
with the familiar vertex operators acting in the Fock space
we introduce the space
$$
  \VGC = \SGC \otimes \mathbb C [ \Rz].
$$

We extend the bilinear form 
$\langle \  , \ \rangle_{\wt }^q$ in $\SGC$ to the space
$\VGC$ and also extend the $\mathbb Z_+$-gradation on $\SGC$
to a $\frac12 \mathbb Z_+$-gradation on $\VG$.

We extend the characteristic map to the map 
\begin{equation*}
ch: \FGC\longrightarrow \VGC
\end{equation*}
by identity on $\Rz$. Then Proposition \ref{P:isometry1}
and Theorem \ref{th_isometry}
imply that we have an isometric isomorphism of Hopf algebras.
We can now identify the operators from the previous
subsections with the operators constructed
from the Heisenberg algebra.

\begin{theorem} 
\label{T:characteristic}
For any $\g \in R(\G)$ and $k\in\mathbb Z$,  we have 
  \begin{eqnarray} \label{E:ch1}
   \ch \bigl ( H_+ (\g\otimes q^k, z) \bigl )
   &=& \exp \biggl ( \sum\limits_{ n \ge 1} \frac 1n \,
    a_{-n} ( \g ) (q^{-k}z)^n \biggr ), \\ \label{E:ch2}
   \ch \bigl ( E_+ (\g\otimes q^k, z) \bigl )
   &=& \exp\biggl ( -\sum\limits_{n\ge 1}\frac 1n  \,
    a_{-n}(\gamma)(q^{-k}z)^n\biggr ),  \\ \label{E:ch3}
   \ch \bigl ( H_- (\g\otimes q^k , z) \bigl )
   &=& \exp \biggl ( \sum\limits_{n \ge 1}\frac 1n \,
    a_n(\g) (q^{-k}z)^{-n}\biggr ),         \\ \label{E:ch4}
   \ch \bigl ( E_- (\g\otimes q^k , z) \bigl )
   &=& \exp\,\, \biggl ( -\sum\limits_{ n \ge 1}
    \frac 1n \,{ a_n( \g )} (q^{-k}z)^{ -n} \biggr ) .
  \end{eqnarray}
\end{theorem}

\begin{proof}
  The first and second identities were essentially
 established in Proposition~\ref{prop_exp} together with
 Lemma~\ref{lem_isom}, where the components are viewed as
operators acting on $\RGC$ or $\SGC$. Note that
$a_n(\g\otimes q^k)=a_n(\g)q^{kn}$.

 We observe from definition that
 the adjoint $*$-action
of $E_+ (\g\otimes q^k, z)$ and $H_- (\g\otimes q^k , z)$
 with respect to the bilinear form
 $\langle \ , \ \rangle_{\wt}^q$
 are $E_- (\g\otimes q^{k} , z^{-1})$ and $ H_- (\g\otimes q^{k} , z^{-1})$
 respectively. 
 The third and fourth identities are obtained by applying
 the adjoint action $*$ to the first two identities.
\end{proof}

\begin{remark} Replacing $\g$ by $-\g$ in (\ref{E:ch1}) and (\ref{E:ch3})
we obtain the equivalent formulas (\ref{E:ch2}) and (\ref{E:ch4}) 
respectively.
\end{remark}

Applying the characteristic map to the
vertex operators $Y^{\pm}(\g, k, z)$, we obtain the 
following group theoretical explanation of
vertex operators acting on the Fock space $\FGC$.

\begin{theorem} \label{T:vertexop}
For any $\g\in\RG$ and $k\in\mathbb Z$, we have 
\begin{align*}     
  &Y^{+}( \g , k, z)\\
  &= \exp \biggl ( \sum\limits_{ n \ge 1} 
  \frac 1n \, \widetilde{a}_{-n} ( \g ) z^n \biggr ) \,
  \exp \biggl ( -\sum\limits_{ n \ge 1}
  \frac 1n \,{ \widetilde{a}_n( \g)} q^{-kn} z^{ -n} \biggr )
  e^{ \g} z^{ \partial_{\g }}\\
&=ch(H_+(\g, z))ch(S(H_+(\g\otimes q^{k}, z^{-1})^*))
e^{ \g} z^{ \partial_{\g }},    
\end{align*}

 \begin{align*}     
  &Y^{-}( \g , k, z)\\
  &= \exp \biggl ( -\sum\limits_{ n \ge 1} 
  \frac 1n \, \widetilde{a}_{-n} ( \g )q^{kn} z^n \biggr ) \,
  \exp \biggl ( \sum\limits_{ n \ge 1}
  \frac 1n \,{ \widetilde{a}_n( \g)}  z^{ -n} \biggr )
  e^{ -\g} z^{-\partial_{\g }}\\
&=ch(S(H_+(\g\otimes q^{k}, z^{-1})))
ch(H_+(\g, z)^*)e^{-\g} z^{-\partial_{\g }}.     
\end{align*}   
\end{theorem}

We note that for $\g\in\G^*, l\in\mathbb Z$
\begin{equation}
Y^{\pm}(\g\otimes q^l, k, z)=Y^{\pm}(\g, k, q^{-l}z).
\end{equation} 

It follows from Theorem \ref{T:vertexop} that
\begin{align*}
&ch\big(Y^{\pm}(\g, k, z)\big)=X^{\pm}(\g, k, z)\\
&= \exp \biggl ( -\sum\limits_{ n \ge 1} 
  \frac 1n \, {a}_{-n} ( \g )q^{n{(k\mp k)}/2} z^n \biggr )\\
 &\qquad\times \exp \biggl ( \sum\limits_{ n \ge 1}
  \frac 1n \,{ {a}_n( \g)}q^{n{(-k\mp k)}/2}  z^{ -n} \biggr )
  e^{ \pm\g} z^{\pm\partial_{\g }}.
\end{align*}

In general the vertex operators
$Y^{\pm}( \g , k, z)$ (for $k\in\mathbb Z$)
generalize the vertex operators considered in \cite{J3} (for $k=\pm 1$). When $q=1$ they specialize to the vertex operators $Y^{\pm}( \g , z)$
 studied in \cite{FJW}.
\section{Basic representations and the McKay 
correspondence}
\label{sect_ade}
\subsection{Quantum toroidal algebras}
Let $Q$ be the root lattice of an affine Lie algebra of simply laced type $A$, $D$, or $E$ with the invariant
form $(\ \ |\ \ )$.
The quantum toroidal algebra $\qta$ is the associative algebra generated by 
$x^{\pm}_{i}(n)$, $a_{i}(m)$, $q^{d}$, $q^c$ ,
$0\leq i\leq r, n, m\in \mathbb Z$ subject to the following relations 
\cite{GKV}:
\begin{gather}
q^{d}a_i(n)q^{-d}=q^na_i(n), 
q^{d}x^{\pm}_i(n)q^{-d}=q^nx_i^{\pm}(n),\\
[a_i(m), a_j(n)]=\delta_{m,-n}\frac{[(\a_i|\a_j)m]}{m}
\frac{q^{mc}-q^{-mc}}{q-q^{-1}},
\label{E:heisenberg}\\ \label{E:comm1}
[a_i(m), x_j^{\pm}(n)] =\pm \frac{[(\a_i|\a_j)m]}{m}q^{\mp |m|c/2}
x_j^{\pm}(m+n),\\ \label{E:comm2}
(z-q^{\pm (\a_i, \a_j)}w)x_i^{\pm}(z)x_j^{\pm}(w)=x_j^{\pm}(w)x_i^{\pm}(z)(q^{\pm (\a_i, \a_j)}z-w), \\ \label{E:comm3}
[x^+_i(z), x^-_j(w)]
=\frac{\delta_{ij}
\{\delta(zw^{-1}q^{-c})\psi_i^+(wq^{c/2})-\delta(zw^{-1}q^c)
\psi_i^-(zq^{c/2})\}}{q-q^{-1}},\\ \label{E:comm4}
Sym_{z_1, \ldots z_N}\sum_{s=0}^{N=1-(\a_i, \a_j)}(-1)^s\bmatrix N\\s\endbmatrix
x^{\pm}_i(z_1)\cdots x_i^{\pm}(z_s)\cdot\\  \nonumber
{\kern .5cm}\cdot x^{\pm}_j(w)x_i^{\pm}(z_{s+1})\cdots
x^{\pm}_i(z_N)=0, \quad\text{for}\quad (\a_i|\a_j)\leq 0, 
\end{gather}
where the generators $\a(n)$ are related to $\psi^{\pm}_{i}(\pm n)$ via:
\begin{gather} \label{E:comm5}
\psi_i^{\pm}(z)
=\sum_{n\geq 0}\psi_i^{\pm}(\pm n)z^{\mp n}=k_i^{\pm 1}exp(\pm(q-q^{-1})
\sum_{n>0}\a_i(\pm n)z^{\mp n}),
\end{gather}
and the Gaussian polynomial 
$$
\bmatrix m\\ n\endbmatrix=\frac{[m]!}{[n]![m-n]!}, \qquad
[n]!=[n][n-1]\cdots [1].
$$
The generating function of $x_n^{\pm}$ are defined by
$$
x_i^{\pm}(z)=\sum_{n\in \mathbb Z}x_i^{\pm}(n)z^{-n-1}, \qquad i=0, \ldots, r.
$$

The quantum toroidal algebra contains a special 
subalgebra-- the quantum affine algebra 
$\qaa$, 
which is generated by simply omitting the
generators associated to $i=0$. The relations are called the Drinfeld
realization of the quantum affine algebras.

In the case of type $A$, the quantum toroidal algebra $\qta$
admits a further deformation $\qtap$.
Let $(b_{ij})$ be the skew-symmetric $(r+1)\times (r+1)$-matrix
\begin{equation}\label{E:skewsymm}
\pmatrix
0 & 1 & 0 & \cdots & 0 & -1\\
-1 & 0 & 1 & \cdots & 0 & 0\\
0 & -1 & 0 & \cdots & 0 & 0\\
\cdots & \cdots & \cdots & \cdots & \cdots & \cdots\\
0 & 0 & 0 & \cdots & 0 & 1\\
1 & 0 & 0 & \cdots & -1 & 0\endpmatrix.
\end{equation}
The quantum toroidal algebra $\qtap$ is the associative algebra generated by 
$x^{\pm}_{in}$, $a_i(m)$, $q^{d_1}$, $q^{d_2}$, $q^c$ ,
$0\leq i\leq r, m, n\in \mathbb Z$ subject to the following relations 
\cite{GKV, VV}:
\begin{gather}
q^{d_1}a_i(n)q^{-d_1}=q^na_i(n), 
q^{d_1}x^{\pm}_i(n)q^{-d_1}=q^nx_i^{\pm}(n),\\
q^{d_2}a_i(n)q^{-d_2}=a_i(n)\\
q^{d_2}x^{\pm}_i(n)q^{-d_2}
=q^{\pm\delta_{n0}}x_i^{\pm}(n),\\
[a_i(m), a_j(n)]=\delta_{m,-n}\frac{[(\a_i|\a_j)m]}{m}
\frac{q^{mc}-q^{-mc}}{q-q^{-1}}p^{mb_{ij}},
\label{E:heisenberg2}\\  \label{E:comm1'}
[a_i(m), x_j^{\pm}(n)] =\pm \frac{[(\a_i|\a_j)m]}{m}q^{\mp |m|c/2}
p^{mb_{ij}}
x_j^{\pm}(m+n),\\  \label{E:comm2'}
(p^{b_{ij}}z-q^{\pm (\a_i|\a_j)}w)x_i^{\pm}(z)x_j^{\pm}(w)
=x_j^{\pm}(w)x_i^{\pm}(z)(p^{b_{ij}}q^{\pm
(\a_i|\a_j)}z-w), \\ \label{E:comm3'}
[x^+_i(z), x^-_j(w)]
=\frac{\delta_{ij}
\{\delta(zw^{-1}q^{-c})\psi_i^+(wq^{c/2})-\delta(zw^{-1}q^c)
\psi_i^-(zq^{c/2})\}}{q-q^{-1}},
\\ \label{E:comm4'}
Sym_{z_1, \ldots z_N}\sum_{s=0}^{N=1-(\a_i|\a_j)}(-1)^s\bmatrix N\\s\endbmatrix
x^{\pm}_i(z_1)\cdots x_i^{\pm}(z_s)\cdot\\ \nonumber
{\kern .5cm}\cdot x^{\pm}_j(w)x_i^{\pm}(z_{s+1})\cdots
x^{\pm}_i(z_N)=0, \quad\text{for}\quad (\a_i|\a_j)\leq 0, 
\nonumber
\end{gather}
where the generators $a_i(n)$ are related to $\psi^{\pm}_{i}(\pm m)$ via:
\begin{gather} \label{E:comm5'}
\psi_i^{\pm}(z)
=\sum_{n\geq 0}\psi_i^{\pm}(\pm n)z^{\mp n}
=k_i^{\pm 1}exp(\pm(q-q^{-1})\sum_{n>0}\a_i(\pm
n)z^{\mp n}).
\end{gather}

We recall that the {\it basic module} of $\qta$ is the 
simple module
generated by the highest weight vector $v_0$ such that
\begin{align*}
&a_i(n+1).v_0=0, 
\qquad x^{\pm}_i(n).v_0=0, \qquad n\geq 0\\ 
& q^c.v_0=qv_0, \qquad q^d.v_0=v_0.
\end{align*}
We say a module is of level one if $q^c$ acts as $q$.

\subsection{A new form of McKay correspondence}

 In this subsection we let $\G$ to be a finite
subgroup of $SU_2$ and consider two distinguished
choices of the class function $\wt$ in $\RGC$ introduced in 
Sect. \ref{S:Mcweights}.

First we consider
$$\wt=\g_0\otimes (q+q^{-1})  - \pi\otimes 1_{\Cs},
$$
where 
$\pi$ is the character of the two-dimensional natural representation 
of $\G$ in $SU_2$.

The Heisenberg algebra in this case has the following relations (cf.
Prop. \ref{prop_orth} and (\ref{E:qcartan1})).
\begin{equation}\label{E:heisen1}
[a_m(\g_i), a_n(\g_j)]=
\begin{cases}
m\delta_{m, -n}(q^m+q^{-m})C, & i=j\\
m\delta_{m, -n}a_{ij}^1C, & i\neq j
\end{cases},
\end{equation}
where $a_{ij}^1$ are the entries of the affine Cartan matrix of ADE type 
(see (\ref{E:qcartan}) at $d=2$).

When $\G\neq \mathbb Z/2\mathbb Z$ or $1$, the relations 
(\ref{E:heisen1}) can be 
simply written as follows:
\begin{equation*}
[a_m(\g_i), a_n(\g_j)]=m\delta_{m, -n}[a_{ij}]_{q^m}C.
\end{equation*}

Recall that the matrix $A^1 = (\langle \g_i, \g_j\rangle_{\wt}^1)
=(a_{ij}^1)_{0 \leq i,j \leq r}$
is the Cartan matrix for the corresponding affine Lie
algebra \cite{Mc}. In particular
$a_{ii}^1 =2$; $a_{ij}^1 =0$ or $-1$ when $i \neq j$ 
and $\G \neq \mathbb Z / 2\mathbb Z$. In the case
of $\G = \mathbb Z / 2\mathbb Z$, $a_{01}^1 =a_{10}^1= -2$.
Let $\mathfrak g$ (resp. $\hat{\mathfrak g}$)
be the corresponding simple Lie
algebra (resp. affine Lie algebra )
associated to the Cartan matrix 
$(a_{ij}^1)_{1\leq i, j\leq r}$ (resp. $A$). Note that the lattice
$\Rz$ is even in this case.

We define the normal ordered product of vertex operators
as follows. 
\begin{align*}
&:Y^+(\g_i, k, z)Y^+(\g_j, k', w):\\
=&H_+(\g_i, z)H(\g_j,w)S(H_+(\g_i\otimes q^k, z^{-1})^*
H_+(\g_j\otimes q^{k'}, w^{-1})^*)\\
&\times e^{\g_i+\g_j}z^{\partial_{\g_i}}
w^{\partial_{\g_j}},\\
&:Y^+(\g_i, k, z)Y^-(\g_j, k', w):\\
=&H_+(\g_i, z)H(-\g_j\otimes q^{-k'}, w)S(H_+(\g_i\otimes q^k, z^{-1})^*
H_+(-\g_j\otimes q^{k'}, w^{-1})^*)\\
&\times e^{\g_i-\g_j}z^{\partial_{\g_i}}
w^{-\partial_{\g_j}}.
\end{align*}
Other normal ordered products are defined similarly.

We introduce
for $a\in \mathbb R$ the following $q$-function: 
\begin{align}\label{E:qfunc}
(1-z)_{q^{2}}^{a}&=\frac{(q^{-a+1}z;q^{2})_{\infty}}
{(q^{a+1}z;q^{2})_{\infty}}=exp\biggl(-\sum_{n=1}^{\infty}\frac{[an]}{n[n]}z^n
\biggl)\\
&=\sum_{m=0}^{\infty}
\begin{bmatrix} a\\ m\end{bmatrix} (-z)^m, \nonumber
\end{align}
where we expand the power series using the $q$-binomial theorem and 
\begin{align*}
\begin{bmatrix} a\\ m\end{bmatrix}&=\frac{(q^a-q^{-a})(q^{a-1}-q^{-a+1})\cdots (q^{a-m+1}-q^{-a+m-1})}{(q^m-q^{-m})(q^{m-1}-q^{-m+1})\cdots (q-q^{-1})},\\
(a; q)_{\infty}&=\prod_{n=0}^{\infty}(1-aq^n).
\end{align*}
When $a$ is a non-negative
integer, $\bmatrix a\\ m\endbmatrix$ equals the Gaussian
polynomial.

The identities in the following theorems are understood as usual
by means of correlation functions (cf. e.g. \cite{FJ, J1}).

\begin{theorem}  \label{th_ope}
Let $\xi=\g_0\otimes (q+q^{-1})-\pi\otimes 1_{\Cs}$.
  Then the  vertex operators 
$Y^{\pm}( \g_i, k, z), Y^{\pm}(-\g_j, k, z)$,  
 $\g_i\in  \G^*, k\in\mathbb Z$ 
acting on the group theoretically defined Fock space
$\FGC$
 satisfy the following relations.
 
\begin{eqnarray*}
 && Y^{\pm}(\g_i, k, z) Y^{\pm}(\g_j, k, w)= \ep (\g_i, \g_j)
  :Y^{\pm}(\g_i, k, z) Y^{\pm}(\g_j, k, w):\\
&&\qquad\times\left\{
  \begin{array}{cc}
1& \mbox{ $\langle \g_i, \g_j\rangle_{\xi}^1=0$}\\
(z-q^{\mp k}w)^{-1}& \mbox{ $\langle \g_i, \g_j\rangle_{\xi}^1=-1$}\\
(z-q^{\mp k-1}w)(z-q^{\mp k+1}w)& 
\mbox{ $\langle \g_i, \g_j\rangle_{\xi}^1=2$}
\end{array},
\right. 
 \end{eqnarray*}
\begin{eqnarray*}
 && Y^{\pm}(\g_i, k, z) Y^{\mp}(\g_j, k,  w) = \ep (\g_i, \g_j)
  :Y^{\pm}(\g_i, k, z) Y^{\mp}(\g_j, k, w):\\
&&\qquad\times\left\{
  \begin{array}{cc}
1& \mbox{ $\langle \g_i, \g_j\rangle_{\xi}^1=0$}\\
(z-w)^{-1}& \mbox{ $\langle \g_i, \g_j\rangle_{\xi}^1=-1$}\\
(z-qw)(z-q^{-1}w)& \mbox{ $\langle \g_i, \g_j\rangle_{\xi}^1=2$}
\end{array},
\right. 
 \end{eqnarray*}

\begin{eqnarray*}
 && Y^{\pm}(\g_i, k, z)Y^{\pm}(-\g_j, -k, w) \\
&&\qquad= \ep (\g_i, \g_j)
  :Y^{\pm}(\g_i, k, z) Y^{\pm}(-\g_j, -k, w):\\
&&\qquad\times\left\{
  \begin{array}{cc}
1& \mbox{ $\langle \g_i, \g_j\rangle_{\xi}^1=0$}\\
(z-q^{\mp k}w) & \mbox{ $\langle \g_i, \g_j\rangle_{\xi}^1=-1$}\\
(z-q^{\mp k-1}w)^{-1}(z-q^{\mp k+1}w)^{-1}
& \mbox{ $\langle \g_i, \g_j\rangle_{\xi}^1=2$}
\end{array},
\right. 
 \end{eqnarray*}

\begin{eqnarray*}
 && Y^{+}(\g_i, k, z)Y^{-}(-\g_j, -k, w)\\
&&\qquad = \ep (\g_i, \g_j)
  :Y^{+}(\g_i, k, z)Y^{-}(-\g_j, -k, w):\\
&&\qquad\times\left\{
  \begin{array}{cc}
1& \mbox{ $\langle \g_i, \g_j\rangle_{\xi}^1=0$}\\
(z-q^{-2k}w)^{-1}& \mbox{ $\langle \g_i, \g_j\rangle_{\xi}^1=-1$}\\
(z-q^{-2k-1}w)(z-q^{-2k+1}w)& \mbox{ $\langle \g_i, \g_j\rangle_{\xi}^1=2$}
\end{array},
\right. 
 \end{eqnarray*}

\begin{eqnarray*}
 && Y^{-}(\g_i, k, z)Y^{+}(-\g_j, -k, w)\\
&&\qquad= \ep (\g_i, \g_j)
  :Y^{-}(\g_i, k, z)Y^{+}(-\g_j, -k, w):\\
&&\qquad\left\{
  \begin{array}{cc}
1& \mbox{ $\langle \g_i, \g_j\rangle_{\xi}^1=0$}\\
(z-w)^{-1}& \mbox{ $\langle \g_i, \g_j\rangle_{\xi}^1=-1$}\\
(z-qw)(z-q^{-1}w)& \mbox{ $\langle \g_i, \g_j\rangle_{\xi}^1=2$}
\end{array}.
\right. 
 \end{eqnarray*}
\end{theorem}

\begin{proof} It is a routine computation to see that:
\begin{eqnarray*}
&& E_- (\g_i\otimes q^k, z) H_+ (\g_j\otimes q^l, w)\\
&=&H_+ (\g_j\otimes q^l, w)E_- (\g_i\otimes q^k, z) 
(1-\frac wz q^{l-k})^{\langle \g_i, \g_j\rangle_{\xi}^1}_{q^2},
\end{eqnarray*}
where the $q$-analog of the power series
$(1-x)_{q^2}^n$ is defined in (\ref{E:qfunc}).

In particular, we have 
\begin{align*}
(1-w/z)_{q^2}&=1-w/z,\\
(1-w/z)^2_{q^2}&=(1-q w/z)(1-q^{-1}w/z).
\end{align*}
Then the theorem is proved by observing that $z^{\g}e^{\partial_{\beta}}
=z^{\langle \g, \beta\rangle_{\xi}^1}e^{\partial_{\beta}}z^{\g}$. 
\end{proof}

\begin{remark} Replacing the vertex operator $Y^{\pm}$ by 
$X^{\pm}$ via the characteristic map $ch$
in the above formulas,
we get the corresponding formulas for vertex operators
$X^{\pm}(\g, k, z)$ acting on 
$\VGC$.
\end{remark}

Now we consider the second distinguished
class function
$$\wt^{q, p}=\g_0\otimes (q+q^{-1})-(\g_1\otimes p+\g_{r}\otimes p^{-1}),$$
when $\G$ is a cyclic group of order $r+1$.

In this case the Heisenberg algebra (\ref{E:heisen}) has 
the following relations according to Prop. \ref{prop_orth} and 
(\ref{E:qcartan2}):
\begin{equation}\label{2ndheisen}
[ a_m(\g_i), {a_n(\g_j)}]=m\delta_{m, -n}
[a_{ij}^1]_{q^m}p^{m b_{ij}}C,
\end{equation}
where $a_{ij}^1$ are the entries of the affine Cartan
matrix of type A and $r\geq 2$.
This is the same 
Heisenberg subalgebra ($c=1$)
in $\qtap$ provided that we identify 
$$a_i(n)=\frac{[n]}na_n(\g_i).
$$

Recall that $(b_{ij})$ is the skew-symmetric
matrix given in (\ref{E:skewsymm}).
We need to slightly modify the definition of the middle term in the vertex operators. For each $i=0, 1, \ldots, r$
we define the modified operator $z^{\partial_{\g, p}}$
on the group algebra $\mathbb C[\Rz]$ by
\begin{equation}
z^{\partial_{\g_i, p}}e^{\beta}
=z^{\langle\g_i, \beta\rangle_{\xi}^1}
p^{-\frac 12\sum_{j=1}^r\langle \g_i, m_j\g_j\rangle_{\xi}^1b_{ij}}e^{\beta},
\end{equation}
where $\beta=\sum_{j}m_j\g_j\in \Rz$.

We then replace the operator $z^{\pm\partial_{\g_i}}$
in the definition of
the vertex operators $Y^{\pm}(\g_i, k, z)$ by the operator
$z^{\pm\partial_{\g_i, p}}$. 
The formulas in Theorems \ref{T:vertexop} remain true
after the term $z^{\pm \partial}$ appearing in
the formulas are modified accordingly. 

The proof of the following theorem is similar to that of
Theorem \ref{th_ope}.

\begin{theorem}  \label{th_ope1}
Let $\G$ be a cyclic group of order $r+1$ and let  
$\xi=\g_0\otimes (q+q^{-1})-(\g_1\otimes p+\g_r\otimes p^{-1})$.
  The  vertex operators 
$Y^{\pm}( \g_i, k, z)$ and $Y^{\pm}(-\g_i, k,  z), 
 \g_i\in  \G^*$ acting on the group theoretically defined Fock space
$\FGC$
 satisfy the following relations.
\begin{eqnarray*}
 && Y^{\pm}(\g_i, k, z) Y^{\pm}(\g_j, k, w)=\ep (\g_i, \g_j)
  :Y^{\pm}(\g_i, k, z) Y^{\pm}(\g_j, k, w):\\
&&\qquad\left\{
  \begin{array}{cc}
1& \mbox{ $\langle \g_i, \g_j\rangle_{\xi}^1=0$}\\
p^{-\frac 12b_{ij}}(z-q^{\mp k}p^{b_{ij}}w)^{-1}& \mbox{ $\langle \g_i, \g_j\rangle_{\xi}^1=-1$}\\
(z-q^{\mp k-1}w)(z-q^{\mp k+1}w)& 
\mbox{ $\langle \g_i, \g_j\rangle_{\xi}^1=2$}
\end{array},
\right. 
 \end{eqnarray*}
\begin{eqnarray*}
&&  Y^{\pm}(\g_i, k, z) Y^{\mp}(\g_j, k,  w)= \ep (\g_i, \g_j)
  :Y^{\pm}(\g_i, k, z) Y^{\mp}(\g_j, k, w):\\
&&\qquad\left\{
  \begin{array}{cc}
1& \mbox{ $\langle \g_i, \g_j\rangle_{\xi}^1=0$}\\
p^{-\frac 12b_{ij}}(z-p^{b_{ij}}w)^{-1}& \mbox{ $\langle \g_i, \g_j\rangle_{\xi}^1=-1$}\\
(z-qw)(z-q^{-1}w)& \mbox{ $\langle \g_i, \g_j\rangle_{\xi}^1=2$}
\end{array},
\right. 
 \end{eqnarray*}

\begin{eqnarray*}
  &&Y^{\pm}(\g_i, k, z)Y^{\pm}(-\g_j, -k, w) \\
&&\qquad =\ep (\g_i, \g_j)
  :Y^{\pm}(\g_i, k, z) Y^{\pm}(-\g_j, -k, w):\\
&&\qquad\left\{
  \begin{array}{cc}
1& \mbox{ $\langle \g_i, \g_j\rangle_{\xi}^1=0$}\\
p^{-\frac 12b_{ij}}(z-q^{\mp k}p^{b_{ij}}w) & \mbox{ $\langle \g_i, \g_j\rangle_{\xi}^1=-1$}\\
(z-q^{\mp k-1}w)^{-1}(z-q^{\mp k+1}w)^{-1}& \mbox{ $\langle \g_i, \g_j\rangle_{\xi}^1=2$}
\end{array},
\right. 
 \end{eqnarray*}

\begin{eqnarray*}
 && Y^{+}(\g_i, k, z)Y^{-}(-\g_j, -k, w)\\
&&\qquad= \ep (\g_i, \g_j)
  :Y^{+}(\g_i, k, z)Y^{-}(-\g_j, -k, w):\\
&&\qquad\left\{
  \begin{array}{cc}
1& \mbox{ $\langle \g_i, \g_j\rangle_{\xi}^1=0$}\\
p^{-\frac 12b_{ij}}(z-q^{-2k}p^{b_{ij}}w)^{-1}& \mbox{ $\langle \g_i, \g_j\rangle_{\xi}^1=-1$}\\
(z-q^{-2k-1}w)(z-q^{-2k+1}w)& \mbox{ $\langle \g_i, \g_j\rangle_{\xi}^1=2$}
\end{array},
\right. 
 \end{eqnarray*}

\begin{eqnarray*}
 && Y^{-}(\g_i, k, z)Y^{+}(-\g_j, -k, w)\\
&&\qquad = \ep (\g_i, \g_j)
  :Y^{-}(\g_i, k, z)Y^{+}(-\g_j, -k, w):\\
&&\qquad\left\{
  \begin{array}{cc}
1& \mbox{ $\langle \g_i, \g_j\rangle_{\xi}^1=0$}\\
p^{-\frac 12b_{ij}}(z-p^{b_{ij}}w)^{-1}& \mbox{ $\langle \g_i, \g_j\rangle_{\xi}^1=-1$}\\
(z-qw)(z-q^{-1}w)& \mbox{ $\langle \g_i, \g_j\rangle_{\xi}^1=2$}
\end{array}.
\right. 
 \end{eqnarray*}
 
\end{theorem}

\begin{remark}
Replacing the vertex operators $Y^{\pm}$ by $X^{\pm}$ via the 
characteristic map $ch$ we obtain the 
corresponding results on the space $\VGC$.
\end{remark}

\subsection{Quantum vertex representations of $\qta$}
For each $i=0$, $\dots$, $r$ let 
\begin{equation*}
\widetilde{a_i}(n)=\frac{[n]}n a_n(\g_i).
\end{equation*}
It follows from (\ref{E:heisen}) and (\ref{E:heisen1}) that
\begin{equation}\label{E:heisenberg1}
[\widetilde{a_i}(m), \widetilde{a_j}(n)]=\delta_{m, -n}\frac{[m\langle \g_i, \g_j\rangle_{\xi}^1]}m [m].
\end{equation}
According to McKay,
the bilinear form $\langle \g_i, \g_j\rangle_{\xi}^1$ is
exactly the same as the invariant form $(\ |\ )$
of the root lattice
of the affine Lie algebra $\loopg$. 
This implies that the
commutation relations (\ref{E:heisenberg1})
 are exactly the commutation relations (\ref{E:heisenberg})
 of the Heisenberg
algebra in $\qta$ if we identify $\widetilde{a_i}(n)$ with ${a_i}(n)$. Thus the Fock space $\SGC$ is
a level one representation for
the Heisenberg subalgebra in $\qta$. 
Under the new variable (by 
identifying $a_i(n)$ with $\widetilde{a_i}(n)$) and after a 
 $q$-shift 
we obtain that
\begin{align*}
&X^{+}(\g_i\otimes q^{-k/2}, k, z)\\
&=\exp \biggl ( \sum\limits_{ n \ge 1} 
  \frac {a_{i} (-n )}{[n]} q^{kn/2}z^n \biggr ) \,
  \exp \biggl ( -\sum\limits_{ n \ge 1}
  \frac {a_i(n)}{[n]} q^{kn/2} z^{ -n} \biggr )
  e^{ \g} z^{ \partial_{\g }} ,\\
&X^{-}(\g_i\otimes q^{-k/2}, k,  z)\\
&=\exp \biggl ( -\sum\limits_{ n \ge 1} 
  \frac {a_{i}(-n )}{[n]}q^{-kn/2} z^n \biggr ) \,
  \exp \biggl ( \sum\limits_{ n \ge 1}
  \frac {a_i(n)}{[n]} q^{-kn/2} z^{ -n} \biggr )
  e^{ -\g} z^{ -\partial_{\g }}.
\end{align*}

The following theorem gives a $q$-deformation
of the new form of McKay correspondence in
\cite{FJW} and provides a direct connection from a finite
subgroup $\G$ of $SU_2$ to the quantum toroidal algebra
$\qta$ of $ADE$ type.
 
\begin{theorem}  \label{T:quantum} Given a finite subgroup $\G$ of $SU_2$,
each of the following correspondence gives a vertex representation
of the quantum toroidal algebra $\qta$ on the Fock space
$\FGC$:
\begin{align*}
x_i^{\pm}(n)& \longrightarrow Y^{\pm}_n(\g_i, -1), \\
a_i(n) &\longrightarrow \frac{[n]}n a_n(\g_i), \qquad q^c\longrightarrow q ;
\end{align*}
or
\begin{align*}
x_i^{\pm}(n)& \longrightarrow Y^{\mp}_n(-\g_i, 1), \\
a_i(n) &\longrightarrow \frac{[n]}n a_n(\g_i), \qquad q^c\longrightarrow q,
\end{align*}
where $i=0, \dots, r$, and $n\in\mathbb Z$.
\end{theorem}
\begin{proof}
Using the usual method of
$q$-vertex operator calculus \cite{FJ, J1} and Theorem \ref{th_ope}
we see that
the vertex operators $Y^{\pm}(\g_i, \pm 1, z)$ satisfy 
relations (\ref{E:comm1}), (\ref{E:comm2}) and (\ref{E:comm4}).
 Observe further  that the above vertex operators at $k=\pm 1$
have the same form as those in the basic representations of the quantum 
affine algebras (see \cite{FJ}). Thus
the relations (\ref{E:comm3}) and (\ref{E:comm5})
are also verified. For each fixed $k= 1$ or $-1$ we have shown that the 
operators
$Y^{\pm}(\g_i, \pm1, z)$ give a level one representation
of the quantum toroidal algebra $\qta$ (see also \cite{Sa, J3}).
\end{proof}

\begin{remark}  Replacing $Y^{\pm}$ by $X^{\pm}$ in the above theorem, 
we obtain a vertex
representation of $\qta$ in the space $\VGC$.
\end{remark}

We can easily get the basic representation of the
quantum affine algebra $\qaa$ on a certain distinguished subspace of
$\FGC$. 

Denote by $\SGGC$ the symmetric algebra generated
 by 
$a_{-n} (\g_i)$, $n >0$, $i =1, \ldots , r$ 
over $\mathbb C[q, q^{-1}]$. 
$\SGGC$ is isometric to $\RGGC$. 

We define
$$
 \FGGC  = \RGGC \otimes \mathbb C [ \Rzz]
 \cong \SGGC\otimes \mathbb C [ \Rzz].
$$
 The space $\VGC$ associated to the lattice
$\Rz$ is isomorphic to the tensor
product of the space $\RGGC$ 
and $\Rz$ as well as the space associated to
the rank $1$ lattice $\mathbb Z \alpha_0$.

\begin{corollary} Given a finite subgroup $\G$ of $SU_2$, each of
the following correspondence gives the basic representation of the
quantum affine algebra $\qaa$ on the Fock space $\FGGC$:
\begin{align*}
x_i^{\pm}(n)& \longrightarrow Y^{\pm}_n(\g_i, -1), \\
a_i(n) &\longrightarrow \frac{[n]}n a_n(\g_i), \qquad q^c\longrightarrow q ;
\end{align*}
or
\begin{align*}
x_i^{\pm}(n)& \longrightarrow Y^{\mp}_n(-\g_i, 1), \\
a_i(n) &\longrightarrow \frac{[n]}n a_n(\g_i), \qquad q^c\longrightarrow q,
\end{align*}
where $i=1, \dots, r$.
\end{corollary}
 
In the case of our second distinguished class function 
$$\xi^{q, p}=\g_0\otimes (q+q^{-1})-
(\g_1\otimes p+\g_r\otimes p^{-1}),
$$
we need to consider the Fock space
\begin{equation*}
 {\tFGGC}  = \RGC \otimes \mathbb C [ \Rz / R^0_{\mathbb Z}]\cong
 \SGC\otimes \mathbb C [\Rzz],
\end{equation*}
where $R^0_{\mathbb Z}$ is the radical of the bilinear form
$\langle \ \ , \ \ \rangle_{\wt}^1$. The correspondence
space for $\tFGGC$
under the characteristic map $ch$ will be denoted $\tVGGC$.

Using similar method as in the proof of Theorem \ref{T:quantum} we derive the
the following theorem.

\begin{theorem} Let $\G$ be a cyclic group of order $r+1\geq 2$ 
and $p=q^{\pm 1}$. 
Each of the following correspondence gives the basic
representation of $\qta$ on $\tFGGC$:
\begin{align*}
x_i^{\pm}(n)& \longrightarrow Y^{\pm}_n(\g_i, -1), \\
a_i(n) &\longrightarrow \frac{[n]}n a_n(\g_i), \qquad q^c\longrightarrow q ;
\end{align*}
or
\begin{align*}
x_i^{\pm}(n)& \longrightarrow Y^{\mp}_n(-\g_i, 1), \\
a_i(n) &\longrightarrow \frac{[n]}n a_n(\g_i), \qquad q^c\longrightarrow q,
\end{align*}
where $i=0, \ldots, r$.
\end{theorem}

\begin{remark} The algebraic picture obtained by
replacing the vertex operator $Y^{\pm}$ by $X^{\pm}$ and $\FGC$ by
$\tVGGC$ in the above Theorem
was given by Sato \cite{Sa}.
\end{remark}

This theorem partly shows 
why the two-para\-meter deformation for $\qta$
is only available in the case of type $A$. It also singles out the
special case of $p=q^{\pm 1}$, where the matrix of the
bilinear form $\langle \ , \ \rangle_{\wt}^{q, q^{\pm 1}}$ is
semi-definite positive (see Sect. \ref{S:Mcweights}) which permits
the factorization of $\FGC$ into $\FGGC$.

We remark that our method can be generalized by replacing 
$R(\G)$ by any finite dimensional Hopf algebra with a Haar measure.
A more general deformation is obtained by replacing $\Cs$
by any torsion-free abelian group. 
In another direction one can replace $\Cs$ by its finite
analog $\mathbb Z/r\mathbb Z$ to study $\qaa$ at 
$r$th roots of unity.

\bigskip

\bibliographystyle{amsalpha}

\end{document}